\documentclass[11pt]{amsart}
\usepackage[a4paper,margin=25mm]{geometry}
\usepackage[T1]{fontenc}
\usepackage{lmodern,microtype}
\usepackage{amsmath,amssymb,amsthm,mathtools}
\usepackage{array,booktabs,longtable,xurl}
\usepackage[bookmarks=false,colorlinks=true,linkcolor=blue,citecolor=blue,urlcolor=blue]{hyperref}
\allowdisplaybreaks[2]
\numberwithin{equation}{section}
\newtheorem{theorem}{Theorem}[section]
\newtheorem{lemma}[theorem]{Lemma}
\newtheorem{proposition}[theorem]{Proposition}

\theoremstyle{definition}

\theoremstyle{remark}

\newtheorem*{ac}{Acknowledgements}
\newcommand{\G}{\mathbb G}
\newcommand{\F}{\mathcal F}
\newcommand{\Z}{\mathbb Z}
\newcommand{\T}{\mathbb T}
\newcommand{\C}{\mathbb C}
\newcommand{\R}{\mathbb R}
\newcommand{\N}{\mathbb N}

\newcommand{\cM}{\mathcal M}
\newcommand{\cN}{\mathcal N}
\newcommand{\cH}{\mathcal H}
\newcommand{\cB}{\mathcal B}
\newcommand{\id}{\operatorname{id}}
\newcommand{\Tr}{\operatorname{Tr}}
\newcommand{\tr}{\operatorname{tr}}

\newcommand{\rank}{\operatorname{rank}}

\newcommand{\Dom}{\operatorname{Dom}}

\newcommand{\norm}[1]{\left\lVert#1\right\rVert}

\title[Hausdorff-Young Inequalities]{Sharp Hausdorff--Young inequalities for some locally compact quantum groups}

\author{Jinsong Wu}
\address{Jinsong Wu, Beijing Institute of Mathematical Sciences and Applications, Beijing, 101408, China}
\email{wjs@bimsa.cn}

\author{Fulin Yang}
\address{Fulin Yang, Beijing Institute of Mathematical Sciences and Applications, Beijing, 101408, China}
\email{fulinyang@bimsa.cn}

\date{}

\begin{document}

\begin{abstract}
In this paper, we show that the usual Hausdorff--Young inequalities for locally compact quantum groups such as the quantum \(E_q(2)\), the real quantum \(az+b\), and the extended quantum \(SU_q(1,1)\), where \(0<q<1\), are sharp.
Moreover, there are no group-like projections with finite Haar weight in these locally compact quantum groups.
We obtain a sharp Hausdorff--Young inequality for Kahng's skew-matrix quantum groups; the optimal constant is strictly less than one for \(1<p<2\); at \(p=2\), Plancherel identity gives equality, and these quantum groups have no nonzero group-like projections with finite Haar weight.
\end{abstract}

\maketitle

\section{Introduction}

Compact quantum groups were first initiated by Woronowicz \cite{Wor87, Wor88} as a quantum extension of compact groups and their duals.
Kustermans and Vaes \cite{KuVaes, KV} constructed the framework of locally compact quantum groups, which extends compact quantum groups, locally compact groups, and their duals.
Euclidean space is an important example of a locally compact group, and its Fourier theory is rich and deep in harmonic analysis.
Babenko \cite{Babenko1961} and Beckner \cite{Beckner1975} proved the sharp Hausdorff--Young inequality for Euclidean spaces, which is profound in harmonic analysis.
Kunze \cite{Kunze1958} and Terp \cite{Terp2017} proved the Hausdorff--Young inequality for locally compact groups.
Klein and Russo \cite{KR} proved the sharp Hausdorff--Young inequality for Heisenberg groups when the conjugate exponent \(p'\) is an even integer.
Cowling, Martini, M\"uller, and Parcet \cite{CMMP2019} studied the Hausdorff--Young inequality for semisimple Lie groups.
Cooney \cite{Cooney} and Caspers \cite{Caspers} proved the Hausdorff--Young inequality for locally compact quantum groups.
Recently, Qiu and Xiong \cite{QX} proved the sharp Hausdorff--Young inequality for noncommutative Euclidean spaces, which is important in the study of the sharp Hausdorff-Young inequalities in the noncommutative setting.

The Hausdorff--Young inequality for locally compact quantum groups \cite{Cooney, Caspers} is sharp for compact quantum groups, since there always exists a group-like projection.
However, this is not true for general locally compact quantum groups.
It is an interesting question to find a locally compact quantum group whose optimal constant for its Hausdorff--Young inequality is strictly less than one.
There are some interesting examples of locally compact quantum groups: real quantum \(E_q(2)\) \cite{DaeleWor, Baaj, Pal}, real quantum \(az+b\) \cite{PuszSoltan}, extended quantum \(SU(1,1)\) \cite{KK}, Kahng's skew-symmetric quantum groups \cite{Kahng}, etc.
We shall investigate the sharp Hausdorff--Young inequalities for these locally compact quantum groups.
We shall show that the usual Hausdorff--Young inequalities are sharp for real quantum \(E_q(2)\) (Theorem \ref{e2:sharp}), real quantum \(az+b\) (Theorem \ref{az:sharp}), and the extended quantum \(SU(1,1)\) (Theorem \ref{su:sharp}).
There are no nonzero group-like projections with finite Haar weight in these locally compact quantum groups.
Inspired by the recent work of Qiu and Xiong \cite{QX}, we show that the optimal constant of the sharp Hausdorff--Young inequalities for Kahng's skew-symmetric quantum groups is strictly less than one (Theorem \ref{thm:maink}).

The paper is organized as follows.
In Section 2, we recall notations in the theory of locally compact quantum groups in the von Neumann algebraic setting and its Fourier theory.
In Sections 3, 4, and 5, we prove the sharpness of the Hausdorff--Young inequality for real quantum \(E_q(2)\), real quantum \(az+b\), and the extended quantum \(SU(1,1)\) respectively.
In Section 6, we obtain the sharp Hausdorff--Young inequality for Kahng's skew-symmetric quantum groups.

\begin{ac}
J. Wu was supported by grants from Beijing Institute of Mathematical Sciences and Applications.
J.~Wu was supported by NSFC (Grant no. 12371124). 
\end{ac}

\section{Preliminaries}

In this section, we shall briefly recall the locally compact quantum groups and their Fourier theory.

\subsection{Locally compact quantum groups}

A locally compact quantum group in the sense of Kustermans and Vaes is a pair $\G=(\cM,\Delta)$, 
where $\cM$ is a von Neumann algebra and
$\Delta:\cM\to \cM\overline\otimes \cM$ is a normal unital $*$-homomorphism satisfying
\[
 (\Delta\otimes\id)\Delta=(\id\otimes\Delta)\Delta.
\]
There are normal semifinite faithful weights $\varphi,\psi$ such that, for $x\ge0$ and $\omega\in \cM_*^+$,
\[
 \varphi((\omega\otimes\id)\Delta(x))=\omega(1)\varphi(x),
 \qquad
 \psi((\id\otimes\omega)\Delta(x))=\omega(1)\psi(x),
\]
where $\cM_*^+$ is the space of all positive normal linear functionals on $\cM$.
The weights are called left and right Haar weights, respectively. 
We write $\cM=L^\infty(\G)$.
The construction and duality results used below are those of \cite{KV}.

Let
\[
 \mathfrak n_\varphi=\{x\in \cM:\varphi(x^*x)<\infty\},
 \qquad
 \mathfrak m_\varphi
 =\operatorname{span}(\mathfrak n_\varphi^*\mathfrak n_\varphi),
\]
and let $(\cH_\varphi,\Lambda_\varphi)$ be the GNS construction with respect to $\varphi$, where $\Lambda_\varphi: \mathfrak n_\varphi \to \cH_\varphi$ is the inclusion.
The canonical left multiplicative unitary $W\in \cM\overline\otimes\cB(\cH_\varphi)$ is specified by
\begin{equation}\label{pre:W}
 W^*(\Lambda_\varphi(a)\otimes\Lambda_\varphi(b))
 =(\Lambda_\varphi\otimes\Lambda_\varphi)
       \bigl(\Delta(b)(a\otimes1)\bigr)
\end{equation}
for $a, b\in \mathfrak n_\varphi$. 
In particular, $\Delta(x)=W^*(1\otimes x)W$ for $x\in \cM$.
The reduced left regular representation is
\begin{equation}\label{pre:lambda}
 \lambda:\cM_*\longrightarrow \cB(\cH_\varphi),\qquad
 \lambda(\omega)=(\omega\otimes\id)(W),
\end{equation}
Let $\widehat \cM$ be the von Neumann algebra generated by $\lambda(\omega)$, $\omega\in \cM_*$.
The dual coproduct is $\widehat\Delta(y)=\Sigma W(y\otimes1)W^*\Sigma$, where $\Sigma$ is the flip of the tensor product.
The predual convolution is $\omega*\eta=(\omega\otimes\eta)\Delta$ and obeys $\lambda(\omega*\eta)=\lambda(\omega)\lambda(\eta)$.
The pair $(\widehat \cM, \widehat\Delta)$ is the dual locally compact quantum group of $(\cM, \Delta)$, denoted by $\widehat \G$.

A locally compact quantum group is unimodular if  $\varphi=\psi$.
A locally compact quantum group is compact if $\varphi=\psi$ is a state.
This is also called compact quantum group.
A compact quantum group is of Kac type if $\varphi$ is tracial.

\subsection{Interpolation spaces}

In this section, we shall recall the Hilsum-Terp $L^p$ space via the Connes' spatial derivatives.
Let $\mathcal T_\varphi$ denote the Tomita algebra.
The elements in $\mathcal T_\varphi$ are entire analytic for $\sigma^\varphi$, and all their complex modular translates belong to $\mathfrak n_\varphi\cap\mathfrak n_\varphi^*$.
Write $\mathcal T_\varphi^2$ for the linear span of its products.
For $x\in \mathcal{T}_\varphi^2$, the symbol $x\varphi$ means the normal functional
\begin{equation}\label{pre:xphi}
 (x\varphi)(y)=\varphi(yx).
\end{equation}
The formula first holds where the product is in the weight domain and then denotes its bounded normal extension. 
It is not an assertion that $\varphi(yx)$ is defined for arbitrary bounded $x,y$.

We use Izumi's parameter $\displaystyle z=-\frac12$, called the left convention in
\cite[Section 2]{Caspers}. 
In this convention the common injections are
\[
 j_1(x)=x\varphi,\qquad j_2(x)=\Lambda_\varphi(x),\qquad
 j_\infty(x)=x.
\]
For precision, the intersection used to match the first two
injections is
\[
 \mathcal I_\varphi
 =\{\omega\in \cM_*:\text{ there is } C_\omega<\infty
       \text{ with }|\omega(y^*)|\le C_\omega \norm{\Lambda_\varphi(y)}
\text{ for all }y\in\mathfrak n_\varphi\}.
\]
For $\omega\in\mathcal I_\varphi$ there is a unique
$\xi_\varphi(\omega)\in \cH_\varphi$ such that
\[
 \omega(y^*)=
 \langle\xi_\varphi(\omega),\Lambda_\varphi(y)\rangle.
\]
In particular $\xi_\varphi(x\varphi)=\Lambda_\varphi(x)$ on
$\mathcal T_\varphi^2$.
The compatible-couple constructions of \cite[Sections 2--3]{Caspers} identify the source and target interpolation spaces, for $\theta=2/r$, $r \geq 2$ and $1/p+1/r=1$, as
\begin{equation}\label{pre:couples}
 (\cM_*,\cH_\varphi)_\theta=L^p(\cM),
 \qquad
 (\widehat \cM,\cH_{\widehat\varphi})_\theta
       =L^r(\widehat \cM).
\end{equation}
The equalities are isometric and use the specified compatible
embeddings, not arbitrary simultaneous realizations of the spaces.

Recall Hilsum's spatial $L^p$ spaces. 
Fix a normal semifinite faithful weight $\phi$ on $\cM'$ and put $d_\varphi=d\varphi/d\phi$. Under the canonical identification with the left interpolation space,
\begin{equation}\label{pre:spatial}
 j_s(x)=\overline{x d_\varphi^{1/s}},\qquad 1\le s<\infty,
 \qquad j_\infty(x)=x.
\end{equation}
Different commutant weights give isometrically equivalent realizations.
The modular automorphisms are $\sigma_t^\varphi(X)=d_\varphi^{it}Xd_{\varphi}^{-it}$.
The scaling group is denoted by $\tau_t$.
We have that 
\begin{equation}\label{eq:covariance}
 \Delta\circ\tau_t
   =(\sigma_t^\varphi\otimes\sigma_{-t}^\psi)\circ\Delta.
\end{equation}

In particular, if $\cM$ has a faithful normal semifinite trace $\tau$ and $\varphi(a)=\tau(D_\varphi^{1/2}aD_\varphi^{1/2})$ for $a\ge0$, we may work in ordinary $L^s(\cM,\tau)$:
\begin{equation}\label{pre:weighted}
 j_s(x)=\overline{xD_\varphi^{1/s}},\qquad
 \norm{j_s(x)}_s^s
 =\tau\!\left(\left|\overline{xD_\varphi^{1/s}}\right|^s\right).
\end{equation}
When $\cM=\cB(\cH)$ these are Schatten norms. 
In general the last expression is not $\varphi(|x|^s)$. 
It is equal to that quantity when $|x|$ commutes with $D_\varphi$, provided the relevant products are defined. 
The analogous fact holds for a general weight when $|x|$ is in its centralizer.

\subsection{Fourier transform}
The canonical dual left Haar weight is normalized by
\begin{equation}\label{pre:dual-GNS}
 \widehat\Lambda(\lambda(\omega))
      =\xi_\varphi(\omega),\qquad \omega\in\mathcal I_\varphi,
\end{equation}
Thus on the space $\mathcal T_\varphi^2$,
\begin{equation}\label{pre:Plancherel-core}
 \widehat\varphi\bigl(\lambda(x\varphi)^*\lambda(x\varphi)\bigr)
       =\varphi(x^*x).
\end{equation}
Let $1\le p\le2$, $p'=p/(p-1)$, with $p'=\infty$ if $p=1$.
The Fourier transform
\begin{align}\label{pre:HY}
 \F_p:L^p(\G)\longrightarrow L^{p'}(\widehat\G)
\end{align}
is defined on the core $\mathcal T^2_\varphi$:
\begin{equation}\label{pre:Fourier-formula}
 \F_p\!\left(\overline{xd_\varphi^{1/p}}\right)
       =\overline{\lambda(x\varphi)\widehat d^{1/p'}}.
\end{equation}
At $p=1$ the right side is $\lambda(x\varphi)$.
At $p=2$ this transform is unitary.
The Hausdorff-Young inequality is that $\|\F_p\|_{p\to p'}\leq 1$.

The noncommutative Hausdorff--Young extension was also developed by Cooney \cite{Cooney}. 
We use Caspers's distinguished left interpolation realization consistently. 

Replacing $\varphi$ by $c\varphi$, $c>0$, replaces its canonically paired dual Haar weight by $c^{-1}\widehat\varphi$.
Under these paired replacements the norm ratio in \eqref{pre:Fourier-formula} is unchanged.

For $1\le p\le2$, define
\[
 C_p(\G)=\sup_{\xi\ne0}
       \frac{\norm{\F_p\xi}_{p'}}{\norm\xi_p}.
\]
A nonzero $\xi$ is an \emph{extremizer} if the quotient is $C_p(\G)$.
A normalized sequence $(\xi_m)_{m \geq 1}$ is \emph{extremizing} if $\norm{\xi_m}_p=1$ and $\norm{\F_p(\xi_m)}_{p'}\to C_p(\G)$.
Sharpness of the inequality and existence of a nonzero extremizer are different assertions.

At $p=2$, all nonzero vectors are extremizers by Plancherel. 
For $1<p<2$, Sections 3,4,5 prove $C_p(\G)=1$ by explicit sequences and Section 6 proves $C_p(\G)<1$.


\subsection{Group-like projections}

Let $P$ be a projection in $L^\infty(\G)$, where $\G$ is a locally compact quantum group.
We say that a nonzero projection $P$ is a group-like projection if $\Delta(P)(1\otimes P)=P\otimes P$. (See \cite[Proposition 4.2]{LWW} for more details).
 If $\varphi(P)<\infty$, then
$\tau_t(P)=\sigma_t^\varphi(P)=P$ for all $t\in\R$.
 If $\tau_t(P)=P$ for all $t\in\R$, then
$\sigma_t^\varphi(P)=P$ for all $t\in\R$, without a Haar-finiteness
assumption. (See \cite[Theorems 3.1, 4.3]{FaaKas} for more details.)

\section{Quantum $E_q(2)$}\label{e2:section}

\subsection{Background}

Throughout this section $0<q<1$ and $\mu=q^{-1}>1$.
Suppose that $(\xi_{j,k})_{j,k\in\Z}$ is an orthonormal basis in $\ell^2(\Z)\otimes\ell^2(\Z)$.
We define operators $v$ and $n$ as follows:
\begin{equation}\label{e2:model}
 v\xi_{j,k}=\xi_{j-1,k},\qquad n\xi_{j,k}=q^j\xi_{j,k+1}.
\end{equation}
The von Neumann algebra generated by $v,n$ is $\cM=\cB(\ell^2(\Z))\overline\otimes L^\infty(\T)$ with a trace $\displaystyle \Tr\otimes \int_{\T} \cdot dz$, where the Haar measure $dz$ of $\T$ has mass $1$. 
Let $\{\delta_j\}_{j\in \Z}$ be an orthonormal basis of $\ell^2(\Z)$.
The trace $\Tr$ on $\cB(\ell^2(\Z))$ can be written as $\displaystyle \Tr(\cdot) =\sum_{j\in \Z}\langle \cdot \delta_j, \delta_j \rangle$.

The comultiplication $\Delta$ on $\cM$ is given by
\begin{equation}
 \Delta(v)=v\otimes v,\quad
 \Delta(n)=n\otimes v\dotplus v^*\otimes n,  
\end{equation}
where $\dotplus$ is the closed normal sum.

For every $0\leq X\in \cM$, the weight $\varphi$ is defined as 
\begin{equation}\label{e2:haar}
 \varphi(X)=\int_\T\sum_{j\in\Z}q^{2j}
     \langle X(z)\delta_j,\delta_j\rangle\,dz.
\end{equation}
Its density $D_\varphi$ affiliated with $\cM$ satisfies $D_\varphi\delta_j=q^{2j}\delta_j$. 
The weight $\varphi$ is both left and right invariant, and is nontracial \cite{Baaj,Pal}.
Then the von Neumann algebra $\cM$ together with comultiplication $\Delta$ and Haar weight $\varphi$ is a locally compact quantum group.

Let $u$ be the polar part of $n$. 
Fourier transformation in the second coordinate identifies $u$ with multiplication by $z$ on $\T$.

The dual locally compact quantum group is $\displaystyle \widehat \cM=\prod_{\kappa\in\Z} \cB(\ell^2(\Z))$.
Let $E^{(\kappa)}_{i,j}$ denote its matrix units. 
Baaj's $\omega^{l,k,j-k}$ extracts row $j$, column $k$ in block
$\kappa=l+2k$. 
The tracial weight $\tau$ and its positive affiliated operator $\widehat d$
\begin{align}\label{eq:duale2}
 \tau=\sum_{\kappa\in\Z}\mu^{-\kappa}\Tr_\kappa,
 \qquad
 \widehat d \delta_j=\mu^{2j-\kappa}\delta_j
       \quad\hbox{in block }\kappa
\end{align}
give $\widehat\varphi=\tau(\widehat d\,\cdot).$
In particular, $ \widehat\varphi(E^{(\kappa)}_{h,h})=\mu^{-2\kappa+2h}.$

\subsection{Sharp Hausdorff-Young inequality}
Let $\{E_{jk}\}_{j,k\in \Z}$ be the system of matrix units of $\cB(\ell^2(\Z))$ such that $E_{jj}$ is a projection onto $\overline{\text{span}}\{\delta_{j}\}$.
For $m\in\Z$, we define
\begin{align}\label{eq:xm}
 x_m = \mathbf1_{\{1\}}(|n|) u^{-m}.
\end{align}
Then
\[
 x_m^*x_m=x_mx_m^*=\mathbf1_{\{1\}}(|n|),\qquad j_s(x_m)=x_m,\qquad
 \norm{x_m}_{s}=1\quad(1\le s\le\infty).
\]
They belong to the square of the Tomita algebra.
The scalar matrix formula preceding Proposition 4.19 of \cite{Baaj}
reads
\[
 \omega^{l,k,j-k}(\lambda(x_m \varphi))
 =(-\mu)^{m}\delta_{l,-m}\delta_{k-j,-m}\overline{A(-j,-m)}.
\]
where
\begin{equation}\label{e2:coeffdef}
 A(m,k)=\int_\T g_\mu(\mu^m z)z^{-k}\,dz,\qquad g_\mu(\zeta)=\prod_{j=0}^{\infty}
       \frac{1+\mu^{-1-2j}\zeta}{1+\mu^{-1-2j}\overline\zeta}.
\end{equation}
The coefficients $A(m,k)$ satisfy that
\begin{enumerate}
\item  $A(m,k)=(-\mu)^k A(m-k,-k)=A(-m,k-m)=(-\mu)^{k-m}A(k,m)$.
\item  $\displaystyle \sum_{k\in\Z}|A(m,k)|^2=1$; this implies that $|A(m,k)|\leq 1$.
\item  $\displaystyle \sum_{m\in\Z}\mu^{2m}|A(m,k)|^2=\mu^{2k}$.
\item  $\displaystyle \lim_{m\to-\infty}A(m,k)=\mathbf1_{\{k=0\}} \quad\hbox{uniformly in }k.$
\end{enumerate}
These are \cite[Lemma A.1 and Propositions A.2--A.3]{Baaj}.

\begin{proposition}\label{prop:dualnorme2}
Suppose that $x_m\in \cM$ is defined in Equation \eqref{eq:xm}.
For $2\le p'<\infty$, we have that 
\begin{equation}\label{e2:exactnorm}
\norm{\F_p(x_m)}_{p'}^{p'}
 =\sum_{k\in\Z}\mu^{-2k}|A(m-k,m)|^{p'}
 =\mu^{-2m}\sum_{k\in\Z}\mu^{2k}|A(k,m)|^{p'}.
\end{equation}
For $p'=\infty$, we have that
\begin{equation}\label{e2:inftynorm}
\norm{\lambda(x_m\varphi)}_\infty
      =\sup_{k\in\Z}|A(k,m)|.
\end{equation}
\end{proposition}
\begin{proof}
By Baaj's characterization of the dual space, we have that 
    \begin{equation}\label{e2:matrixformula}
\lambda(x_m \varphi)\big|_{\kappa=-m+2k}
     =(-\mu)^m A(-k,-m)E^{(-m+2k)}_{k,k-m}
     =A(m-k,m)E^{(-m+2k)}_{k,k-m}.
\end{equation}
By Equation \eqref{eq:duale2}, we have that 
\begin{align*}
    \lambda(x_m \varphi)\widehat d^{1/p'}\big|_{\kappa=-m+2k} =\mu^{-m/p'}A(m-k,m)E^{(-m+2k)}_{k,k-m}.
\end{align*}
This proves the proposition.
\end{proof}

For the following calculation, let \(m\ge0\).
We normalize the Hahn--Exton function as follows:
\[
 J_0(t;q^2)=\sum_{s=0}^\infty
       \frac{(-1)^sq^{s(s+1)}t^{2s}}{(q^2;q^2)_s^2},
 \qquad
 (a;b)_s=\prod_{j=0}^{s-1}(1-a b^j).
\]
Hence
\[
 A(-m,0)=[z^0]\frac{(-q^{m+1}z^{-1};q^2)_\infty}
                         {(-q^{m+1}z;q^2)_\infty}.
\]
The two $q$-binomial expansions, valid absolutely and uniformly on $|z|=1$ because $q^{m+1}<1$, are
\[
 (-tz^{-1};q^2)_\infty
   =\sum_{a\ge0}\frac{q^{a(a-1)}t^az^{-a}}{(q^2;q^2)_a},
 \qquad
 \frac1{(-tz;q^2)_\infty}
   =\sum_{b\ge0}\frac{(-1)^bt^bz^b}{(q^2;q^2)_b}.
\]
Only $a=b$ contributes to the constant coefficient. 
With $t=q^{m+1}$, this gives
\[
 A(-m,0)=\sum_{s\ge0}
  \frac{(-1)^sq^{s(s+1)+2ms}}{(q^2;q^2)_s^2}
 =J_0(q^m;q^2).
\]
This is the Hahn--Exton convention in \cite[Section 2.3]{DeCommerKoelink}.

The same series gives the explicit estimate
\begin{equation}\label{e2:Besselerror}
 |1-J_0(q^m;q^2)|
 \le\frac{q^{2m+2}}
          {(q^2;q^2)_\infty^2(1-q^{2m+2})},
 \qquad
 J_0(q^m;q^2)
 =1-\frac{q^{2m+2}}{(1-q^2)^2}+O_q(q^{4m+6}).
\end{equation}
For the first bound use
$(q^2;q^2)_s\ge(q^2;q^2)_\infty$ and $s(s+1)+2ms\ge(2m+2)s$ for $s\ge1$.
For the remainder in the second expansion, separate $s=1$ and bound the remaining absolutely convergent series by a constant depending only on $q$ times its first power $q^{4m+6}$.
Sharpness therefore holds at a fixed deformation parameter; 
it requires no classical limit.

\begin{theorem}[Sharp Hausdorff-Young Inequality for Quantum $E(2)$]\label{e2:sharp}
Suppose that $\G$ is the quantum $E(2)$.
Then for $1\le p\le2$ and $\displaystyle p'=\frac{p}{p-1}$, we have that
\[
 \norm{\F_p}_{L^p(\G)\to L^{p'}(\widehat\G)}=1,
\]
i.e. the Hausdorff-Young inequality is sharp with optimal constant $1$.
\end{theorem}
\begin{proof}
The upper bound is the general interpolation theorem \cite[Theorems 5.6--5.7]{Caspers}. 
We shall show the sharpness here.
Suppose that $x_m=\mathbf1_{\{1\}}(|n|)u^{-m}$ for $m\in \Z$.
We have that 
\begin{align*}
\left\|x_m D_\varphi^{1/p}\right\|_{p}=\|x_m\|_{p}=1.
\end{align*}
By Proposition \ref{prop:dualnorme2}, we have that 
\begin{align*}
1\geq \norm{\F_p(x_m)}_{p'} \geq  |A(m, m)|=|A(-m, 0)|  =|J_0(q^m;q^2)|
\end{align*}
By Equation \eqref{e2:Besselerror} and letting $m \to\infty$, we see that the theorem holds.
\end{proof}

\subsection{Absence of group-like projections with finite Haar weight}

In this section, we shall show that the locally compact quantum group $\G$ does not have a group-like projection $P$ with $\varphi(P)<\infty$.
\begin{lemma}\label{lem:fixed}
Suppose $\G$ is the quantum $E(2)$.
Then 
\begin{equation}\label{eq:fixed}
 \cM^{\sigma^\varphi}\cap \cM^\tau
   =\ell^\infty(\Z)\overline{\otimes}\mathbb C\mathbf 1.
\end{equation}
Consequently, every scaling-invariant group-like projection has the form
\begin{equation}\label{eq:radial}
 P_I=\sum_{j\in I}E_{jj}\otimes1
     =\mathbf 1_{\{q^j:j\in I\}}(|n|)
 \quad\text{for a nonempty }I\subset\Z.
\end{equation}
The sum is taken in the strong operator topology.
\end{lemma}

\begin{proof}
For $X\in \cM$, let $X_{ij}\in L^\infty(\T)$ be its matrix entries.
By using the fact that $\sigma_t^{\varphi}(X)=D_\varphi^{it} X D_\varphi^{-it}$, we have that 
\[
 (\sigma_t^\varphi(X))_{ij}=q^{2it(i-j)}X_{ij}.
\]
If $X$ is fixed for every $t$, then $X_{ij}=0$ whenever $i\ne j$.
Thus $X$ is diagonal, with entries $f_j\in L^\infty(\T)$.
The scaling action rotates the circle coordinate by $q^{2it}$.
Because $\{q^{2it}:t\in\R\}=\T$, a scaling-fixed $f_j$ is invariant
under all circle rotations and hence is constant almost everywhere.
For example, this follows by taking Fourier coefficients: every
nonzero Fourier mode vanishes.
Conversely, a bounded constant diagonal operator is fixed by both
groups. 
This proves~\eqref{eq:fixed}.

For a projection, each constant diagonal entry is $0$ or $1$.
Therefore every scaling-invariant group-like projection is in the form of \eqref{eq:radial} .
\end{proof}

\begin{proposition}\label{prop:radial}
Let $I\subset\Z$ be nonempty and let $P_I$ be as in~\eqref{eq:radial}.
If $P_I$ is group-like projection, then $I=\Z$, so $P_I=1$.
\end{proposition}

\begin{proof}
Assume that $P_I$ is a group-like projection.
The functional calculus gives
\begin{equation}\label{eq:strong-commute}
\Delta(P_I)=\mathbf 1_{\{q^j:j\in I\}}(|\Delta(n)|).
\end{equation}
Since $\Delta(n)$ is normal, $\Delta(P_I)$ strongly commutes with $\Delta(n)$. 
More precisely,
\begin{equation}\label{eq:domain}
 \Delta(P_I)\Dom(\Delta(n))\subset\Dom(\Delta(n)),\qquad \Delta(n) \Delta(P_I)\zeta=\Delta(P_I)\Delta(n) \zeta
 \quad(\zeta\in\Dom(\Delta(n))).
\end{equation}

Fix $k\in I$. 
For an arbitrary $j\in\Z$, define unit vectors $ \xi=\xi_{j,0}\otimes\xi_{k,0}$ and $\eta=\xi_{j+1,0}\otimes\xi_{k,1}.$
Both have their second tensor factor in the range of $P_I$.
The definition of the group-like projection implies
\begin{equation}\label{eq:eigenvectors}
\Delta(P_I)\xi=\mathbf 1_I(j)\xi,\qquad
\Delta(P_I)\eta=\mathbf 1_I(j+1)\eta,
\end{equation}
where $\mathbf 1_I$ is the characteristic function on $\Z$ for $I$.
Hence
\begin{equation}\label{eq:matrix-computation}
\Delta(n)\xi
 =q^j\xi_{j,1}\otimes\xi_{k-1,0}
  +q^k\xi_{j+1,0}\otimes\xi_{k,1}
 =q^j\xi_{j,1}\otimes\xi_{k-1,0}+q^k\eta.
\end{equation}
The two displayed basis vectors are orthogonal. 
Taking inner products
with $\eta$ in $\Delta(P_I)\Delta(n)\xi=\Delta(n)\Delta(P_I)\xi$, and
\eqref{eq:eigenvectors},
we have that
\begin{align*}
 \langle\eta,\Delta(P_I)\Delta(n)\xi\rangle
   &=\langle \Delta(P_I)\eta,\Delta(n)\xi\rangle
     =\mathbf 1_I(j+1)q^k,\\
 \langle\eta,\Delta(n)\Delta(P_I)\xi\rangle
   &=\mathbf 1_I(j)\langle\eta,\Delta(n)\xi\rangle
     =\mathbf 1_I(j)q^k.
\end{align*}
Since $q^k>0$, we obtain
\[
 \mathbf 1_I(j+1)=\mathbf 1_I(j)\qquad(j\in\Z).
\]
Thus $\mathbf 1_I$ is constant. 
As $k\in I$, its constant value is $1$,
and $I=\Z$.
\end{proof}

\begin{theorem}\label{thm:maine2}
Suppose that $\G$ is the quantum $E(2)$.
Then the scaling-invariant group-like projection in $\cM$ is $1$ and there is no nonzero projection $P\in \cM$ with $\varphi(P)<\infty$ such that
\begin{align*}
 \Delta(P)(1\otimes P)=P\otimes P.
\end{align*}
\end{theorem}
\begin{proof}
If $P$ is scaling invariant and a group-like projection, then it is radial by Lemma~\ref{lem:fixed}. 
By Proposition~\ref{prop:radial}, we see that $P=1$.
If instead $P$ is group-like projection and $\varphi(P)<\infty$, then it is scaling invariant.
Hence $P=1$. 
However, this contradicts the fact
\begin{equation}\label{eq:identity-infinite}
 \varphi(1)=\sum_{j\in\Z}q^{2j}=\infty.
\end{equation}
Therefore the theorem is true.
\end{proof}

\section{The real quantum \(az+b\) group}\label{az:section}

This section concerns the real-parameter family \(0<q<1\) for the quantum \(az+b\) groups and proves the sharpness of the Hausdorff-Young inequality.

\subsection{Background}

Put
\[
 \Gamma_q=\{q^\ell z:\ell\in\Z,\ z\in\T\},\qquad
 \cH=L^2(\Gamma_q)\cong\ell^2(\Z^2).
\]
The radial measure is counting measure and the angular measure is the probability Haar measure. 
The orthonormal basis $\{e_{k, \ell}\}_{k, \ell}$ is defined as follows:
\[
 e_{k,\ell}(q^jz)=\mathbf1_{\{j=\ell\}}z^k.
\]
The affiliated normal generators $a, b$ act on $L^2(\Gamma_q)$ by
\begin{equation}\label{az:generators}
 a e_{k,\ell}=q^k e_{k,\ell-1},\qquad
 b e_{k,\ell}=q^\ell e_{k+1,\ell}.
\end{equation}
They form a regular \(q^2\)-pair in the sense of
\cite[Definition 2.1]{PuszSoltan}; in particular, on the algebraic span
of this basis,
\[
 ab=q^2ba,\qquad ab^*=b^*a.
\]
The comultiplication $\Delta$ is defined by
\begin{equation}\label{az:coproduct}
 \Delta(a)=a\otimes a,\qquad
 \Delta(b)=a\otimes b\dot{+}b\otimes1.
\end{equation}

We have that $\cM=L^\infty(\G)=\cB(\cH)$ with a tracial weight $\Tr$.
Choose the right Haar weight $\psi$ and the corresponding left Haar weight $\varphi$ as
\begin{equation}\label{az:weights}
 \psi(x)=\Tr(|b|x|b|),\qquad
 \varphi(x)=\Tr(D_\varphi^{1/2}xD_\varphi^{1/2}),\qquad
 D_\varphi=q^2|a^{-1}b|^2.
\end{equation}
Here and below, positive products involving unbounded densities are
interpreted as normal semifinite weights. The right Haar formula and
the unitary antipode give the left formula; see
\cite[Section 7.2, Proposition 7.7]{Krajczok}. In these conventions the
unitary antipode sends \(b\) to \(-qa^{-1}b\), so applying it to the
right density produces the left density in \eqref{az:weights}.
Explicitly,
\[
 qa^{-1}b\,e_{k,\ell}=q^{\ell-k}e_{k+1,\ell+1},
 \qquad
 D_\varphi e_{k,\ell}=q^{2(\ell-k)}e_{k,\ell}.
\]
Thus the source Haar weight is nontracial. 

Finite matrix units lie in a common Tomita core. Indeed
\(\sigma_t^\varphi(x)=D_\varphi^{it}xD_\varphi^{-it}\), and every matrix
unit is an entire analytic eigenvector of this action with finite
left and right Haar square norms. The same holds after taking its
adjoint. Consequently the formulas below are valid on a dense
Fourier core and extend by the general Hausdorff--Young theorem.

\subsection{Dual structure}

The self-duality of \(\Gamma_q\) is described by the bicharacter $\chi: \Gamma_q\times \Gamma_q\to \T$ by $\chi(q^kz,q^\ell w)=z^\ell w^k.$

The quantum exponential and its angular Fourier coefficients are
\begin{align}
 F_q(tz)&=\prod_{k=0}^\infty
       \frac{1+q^{2k}tz^{-1}}{1+q^{2k}tz},
       && t\in q^\Z,\quad z\in\T,\label{az:quantum-exponential}\\
 c_n(t)&=\int_\T F_q(tz)z^{-n}\,dz,
 &&F_q(tz)=\sum_{n\in\Z}c_n(t)z^n.\label{az:coefficients}
\end{align}
At removable points the continuous extension is understood. 
The function in \eqref{az:quantum-exponential} is unimodular, and
\begin{equation}\label{az:parseval}
 |c_n(t)|\le1,\qquad \sum_{n\in\Z}|c_n(t)|^2=1.
\end{equation}
Its coefficients are real, as follows by replacing \(z\) with
\(\bar z\) and conjugating.

Let \((\widehat a,\widehat b)\) be the copy of \eqref{az:generators} acting on a second copy of \(\cH\).
The representation theorem of Pusz--So{\l}tan gives the unitary
\begin{equation}\label{az:U}
 U=F_q(\widehat b\otimes b)
   \chi(\widehat a\otimes1,1\otimes a),
 \qquad
 (\id\otimes\Delta)U=U_{12}U_{13}.
\end{equation}
The second leg is the original quantum group. The canonical pairing is
\begin{equation}\label{az:W-orientation}
 W=\Sigma U\Sigma\in \cM\overline\otimes\widehat \cM,
 \qquad
 \lambda(\omega)=(\omega\otimes\id)(W).
\end{equation}
There is no adjoint in \eqref{az:W-orientation}. 
This orientation follows from the universal-dual interpretation of \cite[Theorem 4.2]{PuszSoltan}. 
Krajczok's type I Plancherel description identifies the dual von Neumann algebra with \(B(H)\); the dual quantum group is the opposite quantum \(az+b\) group \cite[Section 7.2, footnote 8]{Krajczok}.
Thus, in the displayed dual generators,
\begin{equation}\label{az:dual-coproduct}
 \widehat\Delta(\widehat a)=\widehat a\otimes\widehat a,\qquad
 \widehat\Delta(\widehat b)
   = \widehat b\otimes\widehat a\dotplus 1\otimes\widehat b.
\end{equation}

\begin{lemma}\label{az:kernel}
The unitary \(W\) acts by
\begin{equation}\label{az:W-kernel}
 W(e_{k,\ell}\otimes e_{r,s})
   =\sum_{n\in\Z}c_n(q^{s+\ell-k-r})
       e_{k+n,\ell-r}\otimes e_{r+n,s-k}.
\end{equation}
\end{lemma}
\begin{proof}
The bicharacter in \eqref{az:U} acts as
\[
 \chi(\widehat a,a)(e_{r,s}\otimes e_{k,\ell})
   =e_{r,s-k}\otimes e_{k,\ell-r}.
\]
For example, the radial part of \(\widehat a\) selects the angular power in the second factor, and its phase translates the second radial coordinate; 
applying the definition of \(\chi\) gives exactly the two shifts displayed. 
The modulus of \(\widehat b\otimes b\) on this vector has eigenvalue \(q^{s+\ell-k-r}\), 
while its phase increments the first index in both tensor factors. 
Expanding the function \(F_q\) by \eqref{az:coefficients} and flipping the legs gives \eqref{az:W-kernel}. 
For the fixed input, the output vectors indexed by \(n\) are orthonormal, and \eqref{az:parseval} proves convergence.
\end{proof}

Putting \(\mu=q^{-1}>1\), and using Baaj's coefficients, we obtain that
\begin{equation}\label{az:Baaj-change}
 g_\mu(\mu^m z)=F_q(q^{1-m}z^{-1}),\qquad
 c_n(q^m)= A(1-m,-n).
\end{equation}
We have that for every \(n\in\Z\),
\begin{equation}\label{az:weighted-sum}
 \sum_{m\in\Z}q^{2m}|c_n(q^m)|^2=q^{2(n+1)}.
\end{equation}

\begin{proposition}\label{az:dual-density}
Under the normalization \eqref{az:weights}, the dual left
Haar weight density is
\begin{equation}\label{az:Dhat}
 \widehat D=q^{-2}|\widehat b|^2,\qquad
       \widehat D e_{r,s}=q^{2(s-1)}e_{r,s}. 
\end{equation}
\end{proposition}
\begin{proof}
Let \(P=|e_{0,0}\rangle\langle e_{0,0}|\). Since
\(PD_\varphi=P\), the functional \(P\varphi\) is the vector state
at \(e_{0,0}\), and \eqref{az:W-kernel} gives
\[
 \lambda(P\varphi)e_{r,s}    =\mathbf1_{\{r=0\}}c_0(q^s)e_{0,s}.
\]
Suppose that $\widehat D=C|\widehat b|^2$.
We have that
\[
 1=\norm{PD_\varphi^{1/2}}_2^2
  =\norm{\lambda(P\varphi)\widehat D^{1/2}}_2^2
  =C\sum_s q^{2s}|c_0(q^s)|^2=Cq^2.
\]
Therefore \(C=q^{-2}\). The dual Haar shape is determined before
this computation by Haar uniqueness; Plancherel determines its scalar.
\end{proof}

\subsection{Sharp Hausdorff-Young Inequality}

For \(m\ge1\), we have that 
\begin{align}
 F_q(q^{1-m}z)&=z^{-m}F_q(q^{m+1}z),\label{az:inversion-product}\\
 A(m, m)=c_{-m}(q^{1-m})&=c_0(q^{m+1})= J_0(q^m;q^2)\longrightarrow1.
 \label{az:alpha-limit}
\end{align}

\begin{theorem}[Sharp Hausdorff-Young inequality]\label{az:sharp}
Suppose that $\G$ is the real quantum $az+b$.
We have that 
\begin{equation}\label{az:sharp-inequality}
 \norm{\F_p}_{L^p(\G)\to L^{p'}(\widehat\G)}=1,
 \qquad 1\le p\le2.
\end{equation}
\end{theorem}
\begin{proof}
The upper estimate is the general interpolation theorem
\cite[Theorems 5.6--5.7]{Caspers}. 
Since \(D_\varphi e_{0,0}=e_{0,0}\), the input is \(x_m D_\varphi^{1/p}=x_m\), a rank-one partial isometry, where $x_m=|e_{m,0}\rangle\langle e_{0,0}|.$
Its norm is one for every Schatten exponent.

The functional \(x_m\varphi\) is \(y\mapsto\langle ye_{m,0},e_{0,0}\rangle\), with inner products linear in the first variable. 
In \eqref{az:W-kernel}, slice the first input \(e_{m,0}\) against the first output \(e_{0,0}\). 
This forces \(r=0\) and \(n=-m\), and gives
\begin{equation}\label{az:witness-Fourier}
\lambda(x_m\varphi)e_{r,s}
 =\mathbf1_{\{r=0\}}c_{-m}(q^{s-m})e_{-m,s-m}.
\end{equation}
The entry from \(e_{0,1}\) to \(e_{-m,1-m}\) equals \(A(m, m)\).
By \eqref{az:Dhat}, the initial basis vector \(e_{0,1}\) has dual Haar density eigenvalue one. 
The same entry of
\(\lambda(x_m\varphi)\widehat D^{1/p'}\) is therefore \(A(m, m)\).
Compression to a row and a column is contractive in every Schatten norm, including the operator norm. 
Hence the Fourier norm is at least \(|A(m, m)|\), which tends to one.
\end{proof}


\subsection{Absence of Group-like projections with finite Haar weight}

In this section, we shall show that the real quantum $az+b$ has no group-like projection $P$ with $\varphi(P)<\infty$.

For each $k\in\Z$, let
\[
 Q_k=\mathbf 1_{\{q^k\}}(|a|)
     =\sum_{\ell\in\Z}|e_{k,\ell}\rangle\langle e_{k,\ell}|.
\]
The projections $Q_k$ are pairwise orthogonal and sum strongly to $1$.
Since $a$ is normal and group-like, the functional calculus applied to~\eqref{az:coproduct} gives
\begin{equation}\label{eq:modulus-coproduct}
 \Delta(|a|)=|a|\otimes|a|.
\end{equation}
The tensor basis diagonalizes the operator on the right, so
\begin{equation}\label{eq:DeltaQ}
 \Delta(Q_0)
 =\mathbf 1_{\{1\}}(|a|\otimes|a|)
 =\sum_{k\in\Z}Q_k\otimes Q_{-k}.
\end{equation}
Consequently
\begin{equation}\label{eq:counterexample0}
 \Delta(Q_0)(1\otimes Q_0)
 =\sum_{k\in\Z}Q_k\otimes Q_{-k}Q_0
 =Q_0\otimes Q_0.
\end{equation}

Equation~\eqref{az:weights} gives
\begin{equation}\label{eq:Q-infinite}
 \varphi(Q_0)=\psi(Q_0)=\sum_{\ell\in\Z}q^{2\ell}=\infty.
\end{equation}

More generally, for every additive subgroup $H\le\Z$, the projection
\begin{equation}\label{eq:subgroup-family}
 P_H=\sum_{k\in H}Q_k
\end{equation}
is group-like projection. 
Indeed, on $Q_k\otimes Q_j$ the required identity is the scalar equality
\[
 \mathbf 1_H(k+j)\mathbf 1_H(j)=\mathbf 1_H(k)\mathbf 1_H(j),
\]
which holds precisely because $H$ is a subgroup. Every such $P_H$
contains $Q_0$, hence has infinite left and right Haar weights.

\begin{lemma}\label{lem:support}
Suppose that $\G$ is the real quantum $az+b$.
Let $P$ be a group-like projection of a locally compact quantum group with $\varphi(P)<\infty$. 
Then
\begin{equation}\label{eq:general-support}
 \sigma_t^\varphi(P)=P,\qquad
 |a|^{2it}P=P=P|a|^{2it}\quad(t\in\mathbb R).
\end{equation}
In particular, $P\le\mathbf 1_{\{1\}}(|a|)$.
\end{lemma}

\begin{proof}
For a projection, $\varphi(P)<\infty$ is equivalent to
$P\in\mathcal N_\varphi$. 
Recall that $\delta=|a|^2,\qquad D_\psi=D_\varphi|a|^2.$
The first assertion is
\cite[Proposition~4.2(3)]{LWW}.
Its conclusion puts $P$ in the centralizer of $\varphi$.
Since it also has finite Haar weight, the functional
\[
 (P\varphi)(X)=\varphi(PXP),\qquad X\in \cM_+,
\]
extends to a bounded normal positive functional of norm $\varphi(P)$.
Thus $P$ belongs to the space $\mathcal L_\varphi$ required in
\cite[Corollary~4.6]{LWW}. That corollary gives the second assertion.
The joint fixed space of the unitaries $|a|^{2it}$ is exactly the spectral subspace of $|a|$ at $1$, proving the last claim.
\end{proof}

\begin{theorem}\label{thm:mainaz}
Suppose that $\G$ is the real quantum $az+b$.
There exists a proper nonzero group-like projection $\mathbf 1_{\{1\}}(|a|)$.
There is no nonzero group-like projection $P$ with $\varphi(P)<\infty$.
\end{theorem}
\begin{proof}
    Suppose $P\ne0$ is group-like projection and $\varphi(P)<\infty$.
By Lemma~\ref{lem:support},
\begin{equation}\label{eq:PunderQ}
 P\le\mathbf 1_{\{1\}}(|a|)=Q_0.
\end{equation}
The modular invariance in~\eqref{eq:general-support} means that
$P$ commutes with all $D_\varphi^{it}$.
On $Q_0\cH$, this density has the simple eigenvalues $q^{2\ell}$ with eigenvectors $e_{0,\ell}$. 
Therefore
\begin{equation}\label{eq:diagonal}
 P=\sum_{\ell\in I}|e_{0,\ell}\rangle\langle e_{0,\ell}|
\end{equation}
for a nonempty subset $I\subset\Z$.
Its finite Haar weight is
\begin{equation}\label{eq:finite-sum}
 \varphi(P)=\sum_{\ell\in I}q^{2\ell}<\infty.
\end{equation}
Since $q^{2\ell}\to\infty$ as $\ell\to-\infty$, this forces $I$
to be bounded below. Let $m=\min I$ and set $r=q^m>0$.
Equation~\eqref{eq:diagonal} implies
\begin{equation}\label{eq:spectral-bound}
 P\le \mathbf 1_{[0,r]}(|b|).
\end{equation}

Note that  $\Delta(P)(P\otimes P)=P\otimes P$. 
Since these are projections, we have that 
\begin{equation}\label{eq:range-bound}
 P\otimes P\le\Delta(P)\le\Delta(\mathbf 1_{[0,r]}(|b|))
   =\mathbf 1_{[0,r]}(|\Delta(b)|).
\end{equation}
The last equality is normal spectral calculus for the affiliated
normal operator $b$ and the normal homomorphism $\Delta$.
If $B=\Delta(b)$, the spectral bound in~\eqref{eq:range-bound}
means that every vector $\zeta$ in the range of $P\otimes P$
belongs to $\Dom(B)$ and satisfies
\begin{equation}\label{eq:norm-upper}
 \|B\zeta\|\le r\|\zeta\|.
\end{equation}

Now take the unit vector $\zeta=e_{0,m}\otimes e_{0,m}$.
It belongs to the range of $P\otimes P$ because $m\in I$.
It also lies in the domains of both algebraic summands of $a\otimes b+b\otimes1$. 
Therefore $B$ acts on it by the algebraic formula, giving
\begin{align}
 B\zeta
 &=(a e_{0,m})\otimes(b e_{0,m})
     +(b e_{0,m})\otimes e_{0,m}\notag\\
 &=q^m\bigl(e_{0,m-1}\otimes e_{1,m}
             +e_{1,m}\otimes e_{0,m}\bigr).\label{eq:computation}
\end{align}
The two displayed basis vectors are orthogonal. Hence
\begin{equation}\label{eq:contradiction}
 \|B\zeta\|^2=2q^{2m}=2r^2>r^2,
\end{equation}
contradicting~\eqref{eq:norm-upper}.
This completes the proof of the theorem.
\end{proof}

\section{Extended quantum \texorpdfstring{$SU_q(1,1)$}{SUq(1,1)}}\label{su:section}

In this section, we shall show that the usual Hausdorff-Young inequality for the extended quantum $SU(1,1)$ is sharp.

\subsection{Background}

Throughout this section, $0<q<1$ is fixed and
\[
 I_q=q^\Z\cup(-q^{\N}),\qquad \N=\{1,2,\ldots\}.
\]
The extended quantum group is the Koelink--Kustermans quantum analogue of the normalizer of $SU(1,1)$ in $SL(2,\C)$ \cite{KK}. 
Its von Neumann algebra is
\begin{equation}\label{su:algebra}
 \cM=L^\infty(\T)\,\overline\otimes\,B(\ell^2(I_q)),
\end{equation}
acting on $\cH_0=L^2(\T)\otimes\ell^2(I_q)$. 
We use normalized Haar measure on $\T$ and counting measure on $I_q$. 
Put
\[
 \mathcal E=\operatorname{span}\{z^m\otimes\delta_t:m\in\Z,\ t\in I_q\}.
\]
On this common algebraic domain the coordinate operators are
\begin{align}
 \alpha_0(z^m\otimes\delta_t)
  &=\sqrt{\operatorname{sgn}(t)+t^{-2}}\,
               z^m\otimes\delta_{qt},\label{su:alpha}\\
 \gamma_0(z^m\otimes\delta_t)
  &=t^{-1}z^{m+1}\otimes\delta_t,\qquad
 e_0(z^m\otimes\delta_t)
  =\operatorname{sgn}(t)z^m\otimes\delta_t.\label{su:gamma}
\end{align}
The convention is $\delta_s=0$ for $s\notin I_q$. Denote their closures
by $\alpha,\gamma,e$. The extension also includes the bounded
self-adjoint partial isometry
\[
 u(z^m\otimes\delta_t)=z^m\otimes\delta_{-t},
\]
with the same zero convention. The operators $\alpha,\gamma,e,u$
generate \eqref{su:algebra} in the sense of affiliation
\cite[Definition 2.3 and Lemma 2.4]{KK}. 
In particular,
\[
 \alpha^*(z^m\otimes\delta_t)
 =\sqrt{\operatorname{sgn}(t)+q^2t^{-2}}\,
                 z^m\otimes\delta_{q^{-1}t};
\]
the coefficient vanishes at $t=-q$. The positive affiliated operator
\begin{equation}\label{su:D}
 D=\gamma^*\gamma,\qquad D\delta_t=t^{-2}\delta_t
\end{equation}
satisfies
$\alpha^*\alpha=D+e$ and $\alpha\alpha^*=q^2D+e$.
The algebraic coproduct formulas are
\begin{align}
 \Delta(\alpha)&\supset
       \alpha_0\odot\alpha_0+q e_0\gamma_0^*\odot\gamma_0,\label{su:coproduct-alpha}\\
 \Delta(\gamma)&\supset
       \gamma_0\odot\alpha_0+e_0\alpha_0^*\odot\gamma_0,
       \qquad \Delta(e)=e\otimes e.\label{su:coproduct-gamma}
\end{align}
Here $\odot$ records the algebraic tensor product on
$\mathcal E\odot\mathcal E$. The containments are intentional:
the closed extensions that define the analytic quantum group are part of
the Koelink--Kustermans construction. In particular, the formal
coproduct of $\gamma^*\gamma$ is not essentially self-adjoint, and
one cannot replace that construction by an unspecified closure.
We give its spectral realization below.

Write $E_{u,v}=|\delta_u\rangle\langle\delta_v|$ and
\[
 \Phi(m,u,v)=z^m\otimes E_{u,v}.
\]
These elements span a $\sigma$-weakly dense $*$-algebra, with
\[
 \Phi(m,u,v)^*=\Phi(-m,v,u),\qquad
 \Phi(m,u,v)\Phi(n,s,t)=\delta_{v,s}\Phi(m+n,u,t).
\]
The left and right Haar weights coincide:
\begin{equation}\label{su:Haar}
 \varphi(X)=\psi(X)
   =\int_\T\Tr\!\left(D^{1/2}X(z)D^{1/2}\right)\,dz
   \qquad(X\ge0).
\end{equation}
Thus the quantum group is unimodular, while its Haar weight is
nontracial. Directly from \eqref{su:D},
\begin{equation}\label{su:modular-source}
 \sigma_t^\varphi(\Phi(m,u,v))
     =|v/u|^{2it}\Phi(m,u,v).
\end{equation}
The scaling group is
\begin{equation}\label{su:scaling}
 \tau_t(\Phi(m,u,v))=q^{2mit}\Phi(m,u,v);
\end{equation}
see \cite[Proposition 4.14]{KK}. It rotates the circle coordinate.

The Haar GNS Hilbert space is
\[
 \mathcal H=L^2(\T)\otimes\ell^2(I_q)\otimes\ell^2(I_q),
 \qquad
 f_{m,u,v}=z^m\otimes\delta_u\otimes\delta_v,
\]
and the GNS representation is $X\mapsto X\otimes1$. The normalization is
\begin{equation}\label{su:GNS}
 \Lambda(\Phi(m,u,v))=|v|^{-1}f_{m,u,v}.
\end{equation}
Indeed $\Phi(m,u,v)^*\Phi(m,u,v)=\Phi(0,v,v)$ has Haar weight
$|v|^{-2}$. The $\Phi(m,u,v)$ are entire for the modular group and
belong to its Tomita algebra. Since
$\Phi(m,u,v)=\Phi(m,u,v)\Phi(0,v,v)$, they also belong to the
squared Tomita core used for the canonical Fourier transform.

\subsection{Explicit Clebsch--Gordan coefficients and the coproduct}

We record the coefficients with their normalization so that every
constant in the sharpness calculation can be traced to a scalar formula.
Let $(a,b;q^2)_k=(a;q^2)_k(b;q^2)_k$. 
Define the regularized series
\begin{equation}\label{su:Psi}
 \Psi(a;b;q^2,z)=
 \sum_{k=0}^{\infty}
 \frac{(a;q^2)_k (bq^{2k};q^2)_\infty}{(q^2;q^2)_k}(-1)^kq^{k(k-1)}z^k.
\end{equation}
For $t\in \pm q^\Z$ and $x\ne0$, put
\[
 \chi(t)=\log_q|t|,\qquad
 \kappa(x)=\operatorname{sgn}(x)x^2,\qquad
 \nu(t)=q^{(\chi(t)-1)(\chi(t)-2)/2},
\]
and define
\[
 \epsilon(x,y)=
 \begin{cases}
 -1,&x>0,\ y<0,\\
 1,&x<0\ \text{or}\ y>0.
 \end{cases}
 \qquad
 c_q=\frac{1}{\sqrt2\,q\,(q^2;q^2)_\infty(-q^2;q^2)_\infty}.
\]
For $t,x,y\in I_q$, set $a_t(x,y)=0$ unless
$\operatorname{sgn}(xy)=\operatorname{sgn}(t)$. On that support, define
\begin{equation}\label{su:a-definition}
\begin{split}
 a_t(x,y)={}&c_q\epsilon(x,y)(-1)^{\chi(t)}
       (-\operatorname{sgn}(y))^{\chi(x)}
       |y|\nu(ty/x)\\
 &{}\times
 \sqrt{\frac{(-\kappa(t),-\kappa(y);q^2)_\infty}
              {(-\kappa(x);q^2)_\infty}}\,
 \Psi\!\left(-\frac{q^2}{\kappa(y)};
                 q^2\kappa(x/y);q^2,q^2\kappa(x/t)\right).
\end{split}
\end{equation}
This is \cite[Definition 3.1]{KK}, with $\epsilon(x,y)$ denoting the sign function called $s(x,y)$ there.
All these coefficients are real.

Here are the precise identities used below. 
For $\theta\in\pm q^\Z$, let
$\ell_\theta=\{(x,y)\in I_q^2:y=\theta x\}$.
The families
\[
 \{a_t|_{\ell_\theta}:
           t\in I_q,\ \operatorname{sgn}(t)=\operatorname{sgn}(\theta)\}
\]
form orthonormal bases of $\ell^2(\ell_\theta)$, and the transposed
coefficient matrices are also unitary
\cite[Propositions 3.2--3.3]{KK}. In particular,
\begin{equation}\label{su:a-orthogonality}
 \sum_{t\in q^\Z}a_t(1,v)^2=1
       \qquad(v\in q^\Z).
\end{equation}
The symmetry formula needed here is
\begin{equation}\label{su:a-symmetry}
 a_t(1,1)=(-1)^{\chi(t)}t^{-1}a_1(1,t)
       \qquad(t\in q^\Z),
\end{equation}
which is a specialization of \cite[Proposition 3.5]{KK}.
Consequently
\begin{equation}\label{su:w}
 w_q(t):=a_t(1,1)^2=t^{-2}a_1(1,t)^2,\qquad
 \sum_{t\in q^\Z}w_q(t)=1,\qquad |a_1(1,v)|\le1.
\end{equation}

For clarity, these same coefficients fix the analytic coproduct.
For $r,s,m\in\Z$ and $t\in I_q$, define functions on
$\T\times I_q\times\T\times I_q$ by
\[
 \eta_{r,s,m,t}(\lambda,x,\mu,y)=
 \begin{cases}
 a_t(x,y)\lambda^{r+\chi(y/t)}\mu^{s-\chi(x/t)},
     &y=\operatorname{sgn}(t)q^m x,\\
 0,&\text{otherwise}.
 \end{cases}
\]
They form an orthonormal basis of $\cH_0\otimes \cH_0$.
The unitary $\mathcal V$ defined by
\[
 \mathcal V:\cH_0\otimes \cH_0
     \longrightarrow L^2(\T)\otimes L^2(\T)\otimes \cH_0,\qquad
 \mathcal V\eta_{r,s,m,t}=z^r\otimes z^s\otimes z^m\otimes\delta_t
\]
gives
\begin{equation}\label{su:analytic-delta}
 \Delta(X)=\mathcal V^*(1\otimes1\otimes X)\mathcal V
       \qquad X\in \cM.
\end{equation}
This is the coproduct of \cite[Definition 3.6 and Proposition 3.8]{KK};
its coassociativity and the Haar invariance are construction theorems
from that reference. Formula \eqref{su:analytic-delta} also specifies the
extensions in \eqref{su:coproduct-alpha}--\eqref{su:coproduct-gamma}.
In particular $\eta_{r,s,m,t}$ is an eigenvector of $\Delta(D)$
with eigenvalue $t^{-2}$.

\begin{lemma}\label{su:a-limits}
On the positive lattice $q^\Z$,
\[
 a_1(1,v)\longrightarrow1\quad(v\to+\infty),\qquad
 a_1(1,v)\longrightarrow0\quad(v\to0).
\]
Moreover $|a_1(1,v)|<1$ at every finite $v\in q^\Z$.
\end{lemma}
\begin{proof}
For the large-argument limit, \eqref{su:a-definition} at $t=x=1$ gives
\[
 a_1(1,v)
 =c_q\,v\nu(v)\sqrt{(-v^2;q^2)_\infty}\,
       \Psi(-q^2/v^2;q^2/v^2;q^2,q^2).
\]
For $v=q^{-K}$, $K\ge0$,
\[
 (-q^{-2K};q^2)_\infty
       =2q^{-K(K+1)}(-q^2;q^2)_K(-q^2;q^2)_\infty.
\]
Using $v\nu(v)=q^{(K^2+K+2)/2}$, we obtain the exact simplification
\begin{equation}\label{su:a-beta}
 a_1(1,q^{-K})=
 \frac{1}{(q^2;q^2)_\infty}
 \sqrt{\frac{(-q^2;q^2)_K}{(-q^2;q^2)_\infty}}\,
      \Psi(-q^{2(K+1)};q^{2(K+1)};q^2,q^2).
\end{equation}
In \eqref{su:Psi}, parameters with modulus at most $1$ have uniformly
bounded numerator products, and $(q^2;q^2)_k\ge(q^2;q^2)_\infty>0$.
The summands at $z=q^2$ are therefore bounded by a constant times $q^{k(k+1)}$. 
Dominated convergence in the series and Euler's product identity give
\[
 \Psi(-q^{2(K+1)};q^{2(K+1)};q^2,q^2)
 \longrightarrow \sum_{k=0}^{\infty}
     \frac{(-1)^k q^{k(k+1)}}{(q^2;q^2)_k}
 =(q^2;q^2)_\infty.
\]
The remaining square-root factor in \eqref{su:a-beta} tends to $1$.

For the small-argument limit, \eqref{su:a-symmetry} gives the simpler
bound
\[
 |a_1(1,v)|=v|a_v(1,1)|\le v,
\]
because \eqref{su:a-orthogonality} at $v=1$ bounds every $|a_v(1,1)|$
by $1$.

It remains to prove strictness at each finite $v$. 
Let $\zeta_v=1\otimes\delta_1\otimes1\otimes\delta_v\in \cH_0\otimes \cH_0$.
The spectral realization \eqref{su:analytic-delta} gives its
$\Delta(D)$ spectral masses as $a_t(1,v)^2$, at eigenvalues $t^{-2}$.
Indeed its expansion in the $\eta$ basis has exactly one pair of
angular labels for each $t$, with coefficient $a_t(1,v)$.
On the algebraic domain,
\[
 \Delta(\gamma)\zeta_v
 =\gamma(1\otimes\delta_1)\otimes\alpha(1\otimes\delta_v)
  +e\alpha^*(1\otimes\delta_1)\otimes\gamma(1\otimes\delta_v).
\]
The two vectors have different angular labels and are orthogonal.
Equations \eqref{su:alpha}--\eqref{su:gamma} yield
\[
 \norm{\Delta(\gamma)\zeta_v}^2
    =(1+v^{-2})+(1+q^2)v^{-2}
    =1+(2+q^2)v^{-2}>1.
\]
If $a_1(1,v)^2=1$, the spectral measure would be concentrated at the
eigenvalue $1$, contradicting this finite quadratic-form expectation.
This proves the strict inequality.
\end{proof}

\subsection{The Fourier transform}

Let $W$ be the canonical multiplicative unitary on
$\mathcal H\otimes\mathcal H$, with
$\Delta(X)=W^*(1\otimes X)W$ in the Haar GNS representation.
For inner products linear in the first variable, write
$\omega_{\xi,\eta}(X)=\langle X\xi,\eta\rangle$ and introduce
\begin{equation}\label{su:Q-definition}
 Q(u,v,n)
 =(\omega_{f_{0,u,1},f_{n,v,1}}\otimes\id)(W^*).
\end{equation}
The matrix entries are explicit:
\begin{equation}\label{su:Q-matrix}
\begin{split}
 \langle Q(u,v,n)f_{a,b,c},f_{\ell,d,s}\rangle
 ={}&\delta_{a-\ell,n}\,
    \delta_{\ell,\chi(ub/(vc))}\,
    \delta_{d,\operatorname{sgn}(bc)svq^a/u}\\
 &{}\times |c/s|\,a_c(u,s)a_b(v,d).
\end{split}
\end{equation}
This is \cite[(4.9), (7.13)]{GKK}. The formula includes the sign
restrictions through the coefficients $a_t$ and through membership
of all radial indices in $I_q$. In particular it determines the
full regular Fourier operator, not a selected irreducible component.

The normal functional associated with a coordinate element is
\[
 \Phi(m,u,v)\varphi
     =|v|^{-2}\,
       \omega_{f_{m,u,1},f_{0,v,1}}\big|_{\cM}.
\]
This follows by evaluating both sides on every $\Phi(n,a,b)$.
Taking the adjoint in \eqref{su:Q-definition} consequently gives
\begin{equation}\label{su:Fourier-core}
 \lambda(\Phi(m,u,v)\varphi)=|v|^{-2}Q(v,u,m)^*.
\end{equation}
Therefore, in the canonical $L^p$ convention,
\begin{equation}\label{su:Fourier-Lp}
 \F_p\!\left(|v|^{-2/p}\Phi(m,u,v)\right)
  =|v|^{-2}\,[Q(v,u,m)^*\widehat d^{1/p'}],
\end{equation}
on this common core. The dual density $\widehat d$ is the spatial
density of the canonically dual Haar weight. 
No independently chosen formal-degree normalization enters \eqref{su:Fourier-Lp}.

For $m\in\Z$, let
\begin{align}
&  P=\Phi(0,1,1),\qquad x_m=\Phi(-m,1,1),\notag\\
& \omega_m=x_m\varphi,\qquad \lambda(\omega_m)=Q(1,1,-m)^*.\label{su:witnesses}
\end{align}
Since $D\delta_1=\delta_1$,
\begin{equation}\label{su:source-norms}
 j_p(x_m)=x_m,\qquad \norm{j_p(x_m)}_p=1
       \quad(1\le p<\infty),\qquad \norm{\omega_m}=1.
\end{equation}
Canonical Plancherel and the slice bound imply
\begin{equation}\label{su:basic}
 \widehat\Lambda(\lambda(\omega_m))=f_{-m,1,1},\qquad
 \widehat\varphi(\lambda(\omega_m)^*\lambda(\omega_m))=1,\qquad \norm{\lambda(\omega_m)}_\infty\le1.
\end{equation}

\begin{lemma}\label{su:modular}
For every $m\in\Z$,
\[
 \widehat\sigma_t(\lambda(\omega_m))=q^{-2mit}\lambda(\omega_m),\qquad
 \widehat{w}_m:=|\lambda(\omega_m)|\in\widehat \cM^{\widehat\varphi}.
\]
For every finite $r\ge2$,
\begin{equation}\label{su:moment-norm}
 \norm{[\lambda(\omega_m)\widehat d^{1/r}]}_r^r
      =\widehat\varphi(\widehat{w}_m^r).
\end{equation}
\end{lemma}
\begin{proof}
The modular element is $\delta=1$, since $\varphi=\psi$.
The canonical dual modular formula
\[
 \widehat\nabla^{it}\Lambda(a)
       =\Lambda(\tau_t(a)\delta^{-it})
\]
is \cite[Section 2, before Lemma 2.7]{KV}.
Equations \eqref{su:scaling} and \eqref{su:basic} show that
$\widehat\nabla^{it}\widehat\Lambda(\lambda(\omega_m))
       =q^{-2mit}\widehat\Lambda(\lambda(\omega_m))$.
Since
$\widehat\Lambda(\widehat\sigma_t(\lambda(\omega_m)))
       =\widehat\nabla^{it}\widehat\Lambda(\lambda(\omega_m))$
and the GNS map is injective, the eigenvector identity follows.
Thus $\lambda(\omega_m)^*\lambda(\omega_m)$ and $\widehat{w}_m$ are fixed by $\widehat\sigma$.
The positive operators $\widehat{w}_m$ and $\widehat d$ strongly commute, giving
\[
 |[\lambda(\omega_m)\widehat d^{1/r}]|^r=\widehat{w}_m^r\widehat d.
\]
Its finiteness follows from $0\le \widehat{w}_m\le1$, $r\ge2$, and $\widehat\varphi(\widehat{w}_m^2)=1$.
\end{proof}

\subsection{Convolution calculation}

\begin{proposition}\label{su:fourth}
For every integer $m$,
\begin{equation}\label{su:R}
\widehat\varphi(\widehat{w}_m^4)
   =\sum_{s\in q^\Z}w_q(s)a_1(1,q^{-m}s)^2.
\end{equation}
In particular,
\[
 0\le \widehat\varphi(\widehat{w}_m^4)<1,\qquad
 \widehat\varphi(\widehat{w}_m^4)\longrightarrow1\ (m\to+\infty),\qquad
 \widehat\varphi(\widehat{w}_m^4)\longrightarrow0\ (m\to-\infty).
\]
\end{proposition}
\begin{proof}
The left ideal property of $\mathcal N_{\widehat\varphi}$ gives
$\lambda(\omega_m)^*\lambda(\omega_m)\in\mathcal N_{\widehat\varphi}$. Hence
\[
 \widehat\varphi(\widehat{w}_m^4)
 =\norm{\widehat\Lambda(\lambda(\omega_m)^*\lambda(\omega_m))}^2
 =\norm{\lambda(\omega_m)^*\widehat\Lambda(\lambda(\omega_m))}^2
 =\norm{Q(1,1,-m)f_{-m,1,1}}^2.
\]
The three Kronecker conditions in \eqref{su:Q-matrix} now force
the output labels to be $(0,q^{-m}s,s)$. Thus
\begin{equation}\label{su:Q-special}
 Q(1,1,-m)f_{-m,1,1}
 =\sum_{s\in q^\Z}s^{-1}a_1(1,s)a_1(1,q^{-m}s)
                       f_{0,q^{-m}s,s}.
\end{equation}
There are no negative-lattice terms, by the support condition
in \eqref{su:a-definition}. The squared coefficients are bounded by
$w_q(s)$, so the series converges in Hilbert norm. Orthogonality and
\eqref{su:w} prove \eqref{su:R}.

At each finite $m$, Lemma \ref{su:a-limits} gives
$a_1(1,q^{-m}s)^2<1$ for every $s$; a probability average of these
numbers is strictly less than $1$.
For $m\to+\infty$ or $m\to-\infty$, the arguments $q^{-m}s$
tend respectively to infinity or zero for each fixed $s$.
Both limits therefore follow from Lemma \ref{su:a-limits} and
dominated convergence under the probability weights $w_q$.
\end{proof}

The same calculation can be expressed as an explicit source convolution.
Define
\begin{equation}\label{su:h}
 h_m=\sum_{s\in q^\Z}
       a_1(1,s)a_1(1,q^{-m}s)\,\Phi(0,q^{-m}s,s).
\end{equation}
The row and column indices of these matrix units are all distinct.
Consequently the series defines a bounded operator with norm at most $1$. 
Its weighted $L^1$ and $L^2$ convergence can be checked without appealing to any formal manipulation of infinite Fourier series.

Indeed, by \eqref{su:a-symmetry},
\begin{equation}\label{su:h-alt}
 a_1(1,s)a_1(1,q^{-m}s)
 =(-1)^m q^{-m}s^2a_s(1,1)a_{q^{-m}s}(1,1).
\end{equation}
Therefore
\begin{align}
 \norm{h_mD}_1
 &=\sum_{s\in q^\Z}s^{-2}
        |a_1(1,s)a_1(1,q^{-m}s)|\notag\\
 &=q^{-m}\sum_{s\in q^\Z}|a_s(1,1)a_{q^{-m}s}(1,1)|
 \le q^{-m},\label{su:h-L1}\\
 \norm{h_mD^{1/2}}_2^2
 &=\sum_{s\in q^\Z}s^{-2}a_1(1,s)^2a_1(1,q^{-m}s)^2
 =\widehat\varphi(\widehat{w}_m^4).\label{su:h-L2}
\end{align}
The first estimate holds by Cauchy-Schwarz inequality, since multiplication by $q^{-m}$ permutes $q^\Z$ and both coefficient-square sums equal $1$.
Finite subsets of $q^\Z$ give truncations converging in both weighted norms. 
Each truncation is in the squared Tomita core.

\begin{proposition}[Convolution Density]\label{su:convolution}
The functional $h_m\varphi$ is normal and
\begin{equation}\label{su:convolution-formula}
 \lambda(h_m\varphi)=\lambda(\omega_m)^*\lambda(\omega_m),\qquad
 \omega_m^\sharp*\omega_m=h_m\varphi,\qquad
 \omega_m^\sharp=(-q)^{-m}\omega_{-m}.
\end{equation}
\end{proposition}
\begin{proof}
For a finite truncation $h_{m,F}$, Plancherel gives $\widehat\Lambda(\lambda(h_{m,F}\varphi))=\Lambda(h_{m,F})$.
The $L^1$ convergence in \eqref{su:h-L1} gives operator-norm convergence of the Fourier transforms to $\lambda(h_m\varphi)$, while Equation \eqref{su:h-L2} gives Hilbert-norm convergence of the GNS vectors.
Closedness of $\widehat\Lambda$ consequently gives
\[
 \widehat\Lambda(\lambda(h_m\varphi))
    =\Lambda(h_m)
    =\sum_{s\in q^\Z}s^{-1}a_1(1,s)a_1(1,q^{-m}s)f_{0,q^{-m}s,s}.
\]
The right-hand side is also $\widehat\Lambda(\lambda(\omega_m)^*\lambda(\omega_m))$ by \eqref{su:Q-special}.
Injectivity of the faithful GNS map proves the first identity.

The antipode formula on the coordinate core in \cite[Proposition 4.13]{KK} gives the last identity in \eqref{su:convolution-formula}. 
It can also be checked by evaluating $\omega_m^\sharp(a)=\overline{\omega_m(S(a)^*)}$ on each
$\Phi(n,u,v)$: 
only $u=v=1$ contributes, and $S(\Phi(n,1,1))=(-q)^n\Phi(n,1,1)$.
Both resulting functionals are normal, and the equality on the antipode core identifies the bounded sharp extension.
Now $\lambda(\omega_m^\sharp*\omega_m)=\lambda(\omega_m)^*\lambda(\omega_m)$.
The left regular representation is injective on $\cM_*$, since slices of $W$ on the second leg generate $\cM$ $\sigma$-weakly. 
This proves the convolution identity.
\end{proof}

This explicit density is useful for checking domains, but the shorter
GNS proof of Proposition \ref{su:fourth} already suffices for sharpness.
In particular, no cancellation of an unbounded Haar density is hidden
in the fourth-moment argument.

\subsection{Sharp Hausdorff-Young inequalities}
In this section, we shall obtain the sharpness of the Hausdorff-Young inequality for the extended quantum $SU_q(1,1)$.
\begin{theorem}\label{su:sharp}
Suppose that $\G$ is the extended quantum $SU_q(1,1)$.
We have that 
\begin{equation}\label{su:sharp-equation}
 \norm{\F_p}_{L^p(\G)\to L^{p'}(\widehat\G)}=1
       \qquad(1\le p\le2).
\end{equation}
\end{theorem}
\begin{proof}
The upper bound is the general canonical interpolation theorem
\cite{Caspers}. 
Recall that 
\[
 0\le \widehat{w}_m\le1,\qquad \widehat\varphi(\widehat{w}_m^2)=1,\qquad
 \widehat\varphi(\widehat{w}_m^4)\longrightarrow1.
\]
For $p'\in[2,\infty)$, Lemma \ref{su:modular} identifies the target
norm to the $p'$th power with $\widehat\varphi(\widehat{w}_m^{p'})$.
If $2\le p'\le4$,
\[
 \widehat\varphi(\widehat{w}_m^4)\le\widehat\varphi(\widehat{w}_m^{p'})\le1.
\]
If $p'\ge4$, consider the probability measure
\[
 \mu_m(E)=\widehat\varphi(\widehat{w}_m^2\mathbf1_E(\widehat{w}_m))
        \qquad(E\subseteq[0,1]\text{ Borel}).
\]
Jensen's inequality for $t\mapsto t^{(p'-2)/2}$ gives
\[
 \widehat\varphi(\widehat{w}_m^{p'})
   =\int(s^2)^{(p'-2)/2}\,d\mu_m(s)
   \ge \widehat\varphi(\widehat{w}_m^4)^{(p'-2)/2}.
\]
Again the upper bound is $1$. 
At $p=1$, or $p'=\infty$,
\[
 \widehat\varphi(\widehat{w}_m^4)\le\norm{\widehat{w}_m}_\infty^2\widehat\varphi(\widehat{w}_m^2)
       =\norm{\lambda(\omega_m)}_\infty^2\le1.
\]
All the desired norm limits follow.

For strictness at $2<p'<\infty$, $\widehat\varphi(\widehat{w}_m^4)<1$ implies
$\widehat\varphi(\widehat{w}_m^2-\widehat{w}_m^4)>0$.
Thus the spectral projection of $\widehat{w}_m$ on $(0,1)$ is nonzero.
On that interval $t^2-t^{p'}>0$, so faithfulness yields
\[
 \widehat\varphi(\widehat{w}_m^{p'})<\widehat\varphi(\widehat{w}_m^2)=1.
\]
This completes the proof of the theorem.
\end{proof}

The proof supplies computable bounds:
\begin{equation}\label{su:explicit-bounds}
 \norm{\F_pj_p(x_m)}_{p'}\ge
 \begin{cases}
 \widehat\varphi(\widehat{w}_m^4)^{1/p'},&2\le p'\le4,\\
 \widehat\varphi(\widehat{w}_m^4)^{(p'-2)/(2p')},&4\le p'<\infty,\\
 \widehat\varphi(\widehat{w}_m^4)^{1/2},&p'=\infty.
 \end{cases}
\end{equation}
At $p=2$, Plancherel gives equality for every nonzero input.
At $p=1$, the theorem establishes sharpness by the explicit sequence.

\subsection{Absence of group-like projections with finite Haar weight}

Let $E_\pm=(1\pm e)/2$. Both are nonzero because both signs occur
in $I_q$. 
By the fact that $e=e^*$, $e^2=1$, we have that 
\begin{equation}\label{eq:counterexample}
 \Delta(E_+)=E_+\otimes E_++E_-\otimes E_-,\qquad
 \Delta(E_+)(1\otimes E_+)=E_+\otimes E_+.
\end{equation}
In particular,
\begin{equation}\label{eq:infinite}
 \varphi(E_+)=\sum_{n\in\Z}q^{-2n}=\infty.
\end{equation}
The projection $E_+$ is also scaling invariant.

\begin{theorem}\label{thm:nogp}
Suppose that $\G$ is the extended quantum $SU_q(1,1)$.
The projection $E_+=(1+e)/2$ is a proper, nonzero, and group-like projection.
No nonzero group-like projection has finite Haar weight.
\end{theorem}
\begin{proof}
Suppose that $P$ is a group-like projection with $\varphi(P)< \infty$.
Then $\tau_t(P)=P=\sigma_t^{\varphi}(P)$.
The scaling group rotates the circle coordinate:
\begin{equation}\label{eq:scaling}
 \tau_t(z^j\otimes E_{p,r})=q^{2jit}z^j\otimes E_{p,r}.
\end{equation}
Since $\{q^{2it}:t\in\mathbb R\}=\T$, scaling invariance makes
$P$ independent of $z$. Modular invariance then makes it block
diagonal on the eigenspaces of $D$.
More explicitly, define
\[
 V_n=\operatorname{span}\{\delta_p:|p|=q^n\}.
\]
For $n\le0$, $V_n$ has dimension $1$; 
for $n\ge1$, it has dimension $2$, with basis $\delta_{q^n},\delta_{-q^n}$.
Thus
\begin{equation}\label{eq:blocks}
 P=1\otimes\bigoplus_{n\in\Z}P_n,\qquad
 \varphi(P)=\sum_{n\in\Z}q^{-2n}\rank(P_n)<\infty,
\end{equation}
where each $P_n$ is an orthogonal projection on $V_n$.
Whenever $P_n\ne0$, one has $q^{-2n}\le \varphi(P)$. 
Consequently $D$ is bounded on $P\cH$, with
\begin{equation}\label{eq:bounded-D}
 DP=PD,\qquad 0\le DP\le \varphi(P)P.
\end{equation}
In particular $\gamma P$, $\alpha P$, and $\alpha^*P$ extend
to bounded operators.
Let
\begin{equation}\label{eq:omega}
 \omega(X)=\varphi(P)^{-1}\varphi(PXP),\qquad X\in \cM.
\end{equation}
This is a normal state with support $P$ and $\omega*\omega=\omega$ which follows from the fact that $P$ is a group-like projection.
On $(P\otimes P)(\cH_0\otimes \cH_0)$, both summands $\gamma\otimes\alpha$ and $e\alpha^*\otimes\gamma$ are bounded, by~\eqref{eq:bounded-D}.
Finite radial truncations by $V_n$ and finite Fourier truncations preserve $P\cH_0$ and give a dense subset of $P\cH_0$ contained in $\mathcal E$. Equation~\eqref{su:coproduct-gamma} holds on tensors from this subset. 
Approximating every vector in $(P\otimes P)(\cH_0\otimes \cH_0)$ by such tensors and using closedness of $\Delta(\gamma)$ proves
\begin{equation}\label{eq:bounded-sum}
 (P\otimes P)(\cH_0\otimes \cH_0)\subset\Dom(\Delta(\gamma)),\qquad
 \Delta(\gamma)(P\otimes P)=\bigl(\gamma\otimes\alpha+e\alpha^*\otimes\gamma\bigr)(P\otimes P).
\end{equation}
In particular the state has a finite expectation of $\Delta(\gamma^*\gamma)=\Delta(D)$.
We have that
\begin{align*}
\omega(D)=& (\omega\otimes\omega)(\Delta(D)) \\
 &=\omega(\gamma^*\gamma)\,\omega(\alpha^*\alpha)
   +\omega(\alpha\alpha^*)\,\omega(\gamma^*\gamma)\\
 &=\omega(D)(\omega(D)+\omega(e))+(q^2\omega(D)+\omega(e))\omega(D),
\end{align*}
i.e.
\begin{align*}
  \omega(D)=(1+q^2)\omega(D)^2+2\omega(D)\omega(e).   
\end{align*}
The omitted cross terms vanish. 
For example the bounded compression $P\gamma^*e\alpha^*P=(\gamma P)^*e\alpha^*P$ transforms under $\tau_t$ by the nonconstant scalar $q^{-2it}$.
Scaling invariance of $\omega$ makes its expectation zero.
The adjoint cross term vanishes as well.
All expressions here are bounded compressions on the support of $\omega$, so no product of uncontrolled unbounded operators is used.

By the fact that $\Delta(e)=e\otimes e$, we have that
\begin{equation}\label{eq:s-idempotent}
 \omega(e)=(\omega\otimes \omega)(\Delta(e))=(\omega\otimes\omega)(e\otimes e)=\omega(e)^2,
\end{equation}
i.e. $\omega(e)=0, 1$.
Note that $ \omega(D)>0$.
We obtain that
\begin{equation}\label{eq:linear-moment}
 1=(1+q^2) \omega(D)+2 \omega(e).
\end{equation}
If $\omega(e)=1$, the right side is strictly greater than $2$, impossible.
If $\omega(e)=0$, the equation forces
\begin{equation}\label{eq:forced-m}
 \omega(D)=\frac1{1+q^2}.
\end{equation}
But on the negative branch $p=-q^n$, $n\ge1$, one has $p^{-2}\ge q^{-2}$. Thus, as positive operators,
\begin{equation}\label{eq:negative-bound}
 D\ge q^{-2}E_-.
\end{equation}
When $\omega(e)=0$, $\omega(E_-)=(1-\omega(e))/2=1/2$, so
\[
  \omega(D)\ge\frac1{2q^2}>\frac1{1+q^2}.
\]
The strict difference is
\[
 \frac1{2q^2}-\frac1{1+q^2}
 =\frac{1-q^2}{2q^2(1+q^2)}>0.
\]
This contradicts~\eqref{eq:forced-m} and proves nonexistence.
\end{proof}

\section{Kahng's skew-matrix quantum groups}

In this section, we shall obtain a sharp Hausdorff-Young inequality for the locally compact quantum groups constructed by Kahng \cite{Kahng}.
We shall apply the idea originating from Qiu and Xiong \cite{QX} and further developed by the authors in \cite{WuYang}.

\subsection{Background}

Let \(n\geq2\), and let \(J\neq0\) be a real skew-symmetric
\(n\times n\) matrix. 
Throughout this section, put
\[
 1<p<2,\qquad p'=\frac p{p-1},\qquad
 B_p=\left(\frac{p^{1/p}}{{p'}^{1/p'}}\right)^{1/2}.
\]
All Euclidean Fourier transforms use the kernel \(e^{-2\pi i x\cdot\xi}\) and the Lebesgue measure. 

Write \(a=(x,y)\in\R^{2n}\), and use \(t\in\R\) for the central
parameter. Set
\begin{align*}
 \sigma_t^J((x,y),(x',y'))
 &=\exp\left[-2\pi i\left(t x\cdot y'
       +\frac{t^2}{2}(Jy)\cdot y'\right)\right],\\
 A_t(x,y)&=\left(x-tJy,y\right),\qquad
 B_t(x,y)=\left(x+\frac t2Jy,y\right).
\end{align*}
Both linear maps $A_t, B_t$ have determinant one. 
Suppose that $ \operatorname*{ess\,sup}_{t}\|f(\cdot,t)\|_1<\infty$.
The representation $L_f$ on $L^2(\R^{2n+1})$ is
\[
 (L_f\xi)(a,t)=\int_{\R^{2n}}
    f(b,t)\sigma_t^J(b,a-b)\xi(a-b,t)\,db.
\]
For a Schwartz function $f\in L^2(\R^{2n+1})$,  the Haar trace and the inclusion map are given by
\begin{equation}\label{eq:source-trace}
 \varphi(L_f)=\int_\R f(0,t)\,dt,
 \qquad \Lambda_{\varphi}(L_f)=f.
\end{equation}
The von Neumann algebra $\cM$ is the direct integral of the twisted group von Neumann algebras \(\cM_t\) with trace $\tr_t$. 
Accordingly, for any $X\in \cM$, we have that 
\[
 \norm X_p^p=\int_\R\tr_t(|X_t|^p)\,dt.
\]
The cocycle identity
\[
 \sigma_t^J(a,b)=\sigma_t^0(B_ta,B_tb)
\]
identifies each fiber with the usual Heisenberg Weyl fiber, preserving its canonical trace. 
The coefficient change is \(f(a,t)\mapsto f(B_t^{-1}a,t)\), followed, if desired, by the unimodular Weyl gauge \(e^{\pi it x\cdot y}\).

The dual von Neumann algebra is the crossed product
\[
 \widehat \cM=L^\infty(\R^{2n})\rtimes\R
\]
for the action induced by \(A_t\). 
For  $Q\in \widehat \cM$, we have that
\begin{equation}\label{eq:Q}
 (Q(h)\,\xi)(a,r)=\int_\R h(a,t)\xi(A_{-t}a,r-t)\,dt.
\end{equation}
Its canonical dual trace satisfies
\begin{equation}\label{eq:dual-trace}
 \widehat\varphi(Q(h))=\int_{\R^{2n}}h(a,0)\,da,
 \qquad \norm{Q(h)}_2=\norm h_{L^2(da\,dt)}.
\end{equation}
The Fourier transform is
\begin{equation}\label{eq:fourier-coeff}
 \F (L_f)=Q(h),\qquad
 h(a,t)=f(-A_{-t}a,-t)
          \sigma_{-t}^J(A_{-t}a,-A_{-t}a).
\end{equation}
Since \((Jy)\cdot y=0\), the last factor is
\(e^{-2\pi it x\cdot y}\). Thus
\begin{equation}\label{eq:mixed-isometry}
 \norm h_{L^p(\R_t;L^q(\R_a^{2n}))}
 =\norm f_{L^p(\R_t;L^q(\R_a^{2n}))}.
\end{equation}

For completeness, Kahng's multiplicative unitary in these coordinates
has the form
\[
 (V\xi)(a,r;b,s)=\sigma_s^J(A_sa,b-A_sa)
                  \xi(A_sa,r+s;b-A_sa,s).
\]
Slicing the second leg by the trace functional associated with \(L_f\)
gives \eqref{eq:fourier-coeff}: at a fixed \(s\), one uses
\[
 \tr_s(U_s(u)L_{f_s})=f(-u,s)\sigma_s^J(u,-u).
\]
The two norms in \eqref{eq:source-trace} and
\eqref{eq:dual-trace} both equal the coefficient \(L^2\)-norm.
There is consequently no additional Plancherel scalar. 

\subsection{A sharp estimate for the real action}
\label{sec:synthesis}

We prove the missing one-dimensional improvement. 
The proof works for an arbitrary measure-preserving action of \(\R\) on a \(\sigma\)-finite measure space \((S,\mu)\); only the shear action is needed here. 
Let \(\cN=L^\infty(S)\).
Let \(\cN\rtimes\R\) have its canonical semifinite trace, and write
\[
 Q(h)=\int_\R M_{h_t}u_t\,dt,
\]
where $h_t(\cdot)=h(\cdot, t)$ and $u_t$ is the unitary operator implementing the crossed product.
Plancherel and the operator-norm endpoint give
\begin{align}
 \norm{Q(h)}_2=\norm h_{L^2(\R;L^2(S))},\label{eq:endpoint2}\quad
 \norm{Q(h)}_\infty \leq\norm h_{L^1(\R;L^\infty(S))}.
\end{align}
Complex interpolation therefore gives the initial estimate
\begin{equation}\label{eq:rough}
 \norm{Q(h)}_{p'}\leq\norm h_{L^p(\R;L^{p'}(S))}.
\end{equation}

\begin{lemma}\label{lem:synthesis}
For \(1<p<2\),
\begin{equation}\label{eq:synthesis}
 \norm{Q(h)}_{p'}\leq B_p\norm h_{L^p(\R;L^{p'}(S))}.
\end{equation}
\end{lemma}

\begin{proof}
We separate the argument into localization, the frozen coefficient,
and the approximation error.

\emph{Localization.} For a semifinite tracial algebra, let
\(\Omega(Y)=\norm Y_{p'}^p\), \(Y\neq0\). The trace-Hessian
inequality used in \cite[Lemma 3.4]{QX} is
\begin{equation}\label{eq:hessian}
 D^2\Omega(Y)[Z,Z]\geq(p-1)D^2\Omega(Y)[iZ,iZ],
\end{equation}
where $\displaystyle D^2\Omega(Y)[Z,Z]=\left.\frac{\partial^2}{\partial s\partial t}\Omega(Y+sZ+tZ)\right|_{s=t=0}$.
Initially take a smooth compactly supported coefficient \(h\).
Put
\[
 Y(a,b)=Q(h_t e^{-at^2+bt}),\quad
 V(a,b)=\norm{Y(a,b)}_{p'}^p,\quad a\geq0,\ b\in\R.
\]
Multiplication of coefficients by \(e^{ivt}\) is the dual action on the crossed product. 
It preserves the trace and every \(L^{p'}\)-norm.
Writing
\(Z=Q(t h_t e^{-at^2+bt})\) and
\(C=Q(t^2h_t e^{-at^2+bt})\), differentiation gives
\[
 D^2\Omega(Y)[iZ,iZ]=D\Omega(Y)[C],
\]
where $\displaystyle D\Omega(Y)[C]=\left.\frac{\partial}{ \partial t}\Omega(Y+tC)\right|_{t=0}$.
The chain rule and \eqref{eq:hessian} now imply
\begin{equation}\label{eq:heat-subsolution}
 (\partial_b^2+p\partial_a)V(a, b)\geq0.
\end{equation}
Integrating this inequality against the positive heat kernel yields,
for every \(A>0\),
\begin{equation}\label{eq:localization}
 \norm{Q(h)}_{p'}^p\leq
 \left(\frac{pA}{\pi}\right)^{1/2}
 \int_\R\norm{Q(h_t e^{-A(t-s)^2})}_{p'}^p\,ds.
\end{equation}
Indeed, if \(\displaystyle (H_vg)(x)=\left(\frac{1}{4\pi v}\right)^{1/2} \int_{\R} e^{-|x-y|^2/(4v)} g(y) dy \) is the ordinary heat semigroup, the derivative of \(H_{a/p}[V(a,\cdot)](0)\) is nonnegative by
\eqref{eq:heat-subsolution}; evaluating at \(a=0,A\) and substituting
\(b=2As\) gives \eqref{eq:localization}. Compact support in \(t\)
bounds all relevant derivatives by an exponential in \(|b|\),
justifying the heat integrations. The coefficient is nonzero exactly
when its \(L^2\) synthesis is nonzero, so no differentiability issue
at \(Y=0\) occurs on this core.

\emph{Frozen coefficients.} 
We need the product estimate
\begin{equation}\label{eq:kss}
 \norm{M_b\lambda(g)}_{p'}
 \leq\norm b_{L^{p'}(S)}\norm{\widehat g}_{L^{p'}(\R)},
 \qquad {p'}\geq2.
\end{equation}
Here \(\displaystyle \lambda(g)=\int g(t)u_t\,dt\). 
To check the constant without a trace-convention ambiguity, consider the bilinear map
\((b,m)\mapsto M_b m(D)\), where \(m=\widehat g\).
At exponent two, \eqref{eq:endpoint2} gives exactly
\(\norm{M_bm(D)}_2=\norm b_2\norm m_2\); at infinity the
operator product gives the same inequality with constant one.
Bilinear complex interpolation proves \eqref{eq:kss}.

Let \(g_A(t)=e^{-At^2}\) and
\(c_A=(pA/\pi)^{1/(2p)}\). Since translating \(g_A\) contributes
a right unitary \(u_s\), \eqref{eq:kss} gives
\[
 c_A\norm{Q(h_s g_A(t-s))}_{p'}
 \leq c_A\norm{\widehat g_A}_{p'}\norm{h_s}_{p'}
 =B_p\norm{h_s}_{p'}.
\]
The final equality follows directly from
\(\widehat g_A(\xi)=(\pi/A)^{1/2}e^{-\pi^2\xi^2/A}\).

\emph{Approximation error.} Applying \eqref{eq:rough} to
\((h_t-h_s)g_A(t-s)\) and integrating in \(s\), we obtain
\begin{align*}
 &\int_\R c_A^p
       \norm{Q((h_t-h_s)g_A(t-s))}_{p'}^p\,ds\\
 &\qquad\leq
 \left(\frac{pA}{\pi}\right)^{1/2}
 \int_\R e^{-pAv^2}
    \norm{h(\cdot+v)-h}_{L^p(\R;L^{p'}(S))}^p\,dv
 \longrightarrow0.
\end{align*}
This is the approximate-identity property and strong continuity of
translations in the Bochner \(L^p\)-space. 
Minkowski's inequality in the \(s\)-variable, together with \eqref{eq:localization}, now gives
\eqref{eq:synthesis}. 
Finally \eqref{eq:rough} and density extend the conclusion to every coefficient in the stated mixed space.
\end{proof}

\begin{proposition}\label{prop:hy}
For every \(1<p<2\), we have that 
\begin{equation}\label{eq:main}
 \norm{\F_{p}(X)}_{L^{p'}(\widehat\G)}
 \leq B_p^{2n+1}\norm X_{L^p(\G)}.
\end{equation}    
\end{proposition}
\begin{proof}
The fiberwise inverse Weyl inequality, in canonical trace
normalization, is
\begin{equation}\label{eq:fiberHY}
 \norm{f(\cdot,t)}_{L^{p'}(\R^{2n})}
 \leq B_p^{2n}\norm{X_t}_{L^p(M_t,\tau_t)}.
\end{equation}
This follows from \cite[Theorem 4.4]{QX}, because the determinant-one
map \(B_t\) and the Weyl gauge preserve the coefficient norms.
Integration in \(t\) yields
\begin{equation}\label{eq:integrated-fiber}
 \norm f_{L^p_tL^{p'}_a}\leq B_p^{2n}\norm X_p.
\end{equation}
Combining Lemma \ref{lem:synthesis}, \eqref{eq:mixed-isometry}, and \eqref{eq:integrated-fiber}, we get
\[
 \norm{\F (X)}_{p'}
 =\norm{Q(h)}_{p'}
 \leq B_p\norm h_{L^p_tL^{p'}_a}
 =B_p\norm f_{L^p_tL^{p'}_a}
 \leq B_p^{2n+1}\norm X_p.
\]
The coefficient-core identity extends by density and the displayed
bounded maps. 
This proves the proposition.
\end{proof}

\subsection{Sharpness of the Hausdorff-Young inequality}
\label{sec:gaussian}
In this section, we shall prove the sharpness of the Hausdorff-Young inequality obtained in Proposition \ref{prop:hy}.
\begin{theorem}\label{thm:maink}
For every \(1<p<2\), we have that 
\[
 \norm{\F_{p}}_{L^p(\G)\to L^{p'}(\widehat\G)}=B_p^{2n+1}.
\]
\end{theorem}
\begin{proof}
For \(\varepsilon>0\), define
\begin{equation}\label{eq:gaussian}
 f_\varepsilon(x,y,t)=
 e^{-\pi itx\cdot y}
 \exp\left[-\pi\left(
       \left|x+\frac t2Jy\right|^2+|y|^2+
       \frac{t^2}{\varepsilon^2}\right)\right],
 \qquad X_\varepsilon=L_{f_\varepsilon}.
\end{equation}
These coefficients belong to the relevant \(L^p\) domains; the
following exact norm formulas also verify this directly. Their
transform has the particularly convenient coefficient
\begin{equation}\label{eq:fixedcoeff}
 \F (X_\varepsilon)=Qf_\varepsilon.
\end{equation}
In fact,
\(B_{-t}(-A_{-t}a)=-B_ta\), and the gauge factor in
\eqref{eq:fourier-coeff} combines with that in
\eqref{eq:gaussian} to give the same coefficient.

For \(r>1\), set
\begin{equation}\label{eq:R}
 R_r(s)=\frac{s}{(1+s/2)^{r}-|1-s/2|^{r}}\quad(s>0),
 \qquad R_r(0)=\frac1r.
\end{equation}

After the shear and the gauge have been removed, the fiber is
\(e^{-\pi(t/\varepsilon)^2}W_t(e^{-\pi|a|^2})\).
In the Schr\"odinger realization,
\[
 (U_t(x,y)\psi)(u)=e^{-2\pi it y\cdot u}\psi(u+x),
 \qquad W_t(x,y)=e^{-\pi itx\cdot y}U_t(x,y).
\]
The Gaussian operator has integral kernel
\[
 K_t(u,v)=\exp\left[-\pi|u-v|^2
                     -\frac{\pi t^2}{4}|u+v|^2\right].
\]
The one-dimensional Mehler formula gives its eigenvalues, and tensor
products give the \(n\)-dimensional eigenvalues:
\[
 (1+|t|/2)^{-n}
 \left(\frac{1-|t|/2}{1+|t|/2}\right)^{k_1+\cdots+k_n},
 \qquad k_j=0,1,2,\ldots.
\]
Here \(\tau_t=|t|^n\Tr\). Taking absolute values when the ratio
is negative and summing the geometric series therefore gives
\begin{equation}\label{eq:sourceexact}
 \norm{X_\varepsilon}_p^p
 =\int_\R e^{-\pi p(t/\varepsilon)^2}R_p(|t|)^n\,dt.
\end{equation}
Since \(R_p(s)\leq C_p(1+s)^{2-p}\) for \(1<p<2\), the
substitution \(t=\varepsilon v\) and dominated convergence show
\begin{equation}\label{eq:sourceasym}
 \norm{X_\varepsilon}_p^p
 \sim\varepsilon p^{-n-1/2},\qquad
 \norm{X_\varepsilon}_p
 \sim\varepsilon^{1/p}p^{-(2n+1)/(2p)}.
\end{equation}

For almost every \(y\), put \(c=|Jy|>0\), \(e=Jy/c\), and
write \(x=z+se\), with \(z\perp e\). The exceptional set
\(\ker J\) has Lebesgue measure zero because \(J\neq0\).
The shear orbits are the \(s\)-lines. Up to a multiplication
unitary, the fiber of \(Qf_\varepsilon\) on such an orbit is
\(e^{-\pi(|z|^2+|y|^2)}\) times the integral operator with kernel
\begin{equation}\label{eq:targetkernel}
 \frac1c\exp\left[-\pi\left(
             \frac{(s+s')^2}{4}
             +\frac{(s'-s)^2}{c^2\varepsilon^2}\right)\right].
\end{equation}
The phase disappears by conjugation because \(e\cdot y=0\).
The crossed-product trace is \(c\,dy\,dz\,\Tr\) in these
coordinates: the diagonal of the orbit kernel is \(h(x,y,0)/c\),
which recovers \eqref{eq:dual-trace}.

The eigenvalues of \eqref{eq:targetkernel} are
\[
 \frac{\varepsilon}{1+c\varepsilon/2}
 \left(\frac{1-c\varepsilon/2}{1+c\varepsilon/2}\right)^k,
 \qquad k=0,1,2,\ldots.
\]
Summing their absolute \({p'}\)-powers and integrating \(z\) gives
the exact expression
\begin{equation}\label{eq:targetexact}
 \norm{\F (X_\varepsilon)}_{p'}^{p'}
 =\varepsilon^{p'-1}{p'}^{-(n-1)/2}
     \int_{\R^n}e^{-\pi p'|y|^2}R_{p'}(\varepsilon|Jy|)\,dy.
\end{equation}
For \({p'}\geq2\), convexity gives \(0<R_{p'}(s)\leq1/p'\).
Thus dominated convergence gives
\begin{equation}\label{eq:targetasym}
 \norm{\F( X_\varepsilon)}_{p'}^{p'}
 \sim\varepsilon^{p'-1}{p'}^{-n-1/2},\qquad
 \norm{\F ( X_\varepsilon)}_{p'}
 \sim\varepsilon^{1/p}{p'}^{-(2n+1)/(2{p'})}.
\end{equation}
Combining \eqref{eq:sourceasym} and \eqref{eq:targetasym},
\begin{equation}\label{eq:ratio}
 \lim_{\varepsilon\downarrow0}
 \frac{\norm{\F( X_\varepsilon)}_{p'}}{\norm{X_\varepsilon}_p}
 =\frac{p^{(2n+1)/(2p)}}{{p'}^{(2n+1)/(2p')}}=B_p^{2n+1}.
\end{equation}
This proves sharpness for every \(1<p<2\). 

\end{proof}

A nonzero group-like projection with finite Haar weight would force the Hausdorff–Young constant to be \(1\), by \cite[Proposition 4.9]{LWW}, contradicting \(B_p^{2n+1}<1\).

\subsection{Extremizers}
In this section, we shall prove there are no nontrivial extremizers for the sharp Hausdorff-Young inequality.

\begin{lemma}\label{lem:strict}
Let \(r>2\), \(p=r'\), and let \(K\) be a self-adjoint (unbounded) operator
without eigenvalues in $\cH$. 
Suppose
\(0\neq Y\in\mathcal{S}^r\cap C^1(K)\) and
\(Z=[K,Y]\in\mathcal{S}^r\), 
where 
\begin{align*}
C^1(K)=\{ X\in \cB(\cH): \text{ the map } t\mapsto e^{itK}Xe^{-itK} \text{ is differentiable}\}.
\end{align*}
Then, for \(\Omega(Y)=\norm Y_{r}^p\),
\[
 D^2\Omega(Y)[Z,Z]>(p-1)D^2\Omega(Y)[iZ,iZ].
\]
\end{lemma}

\begin{proof}
We analyze the equality case of the trace-Hessian proof in
\cite[Lemmas 3.2--3.4]{QX}. 
Write \(\Phi(Y)=r^{-1}\Tr(|Y|^r)\) and \( \displaystyle \alpha=\frac{1}{\Tr(|Y|^r)}D\Phi(Y)[Z]\), which is real.
Put \(\displaystyle \beta=\frac{1}{\Tr(|Y|^r)} D\Phi(Y)[iZ]\),
 and
\begin{align*}
\Delta_\Phi(Z-\alpha Y)=D^2\Phi(Y)[Z-\alpha Y,Z-\alpha Y]-(p-1)D^2\Phi(Y)[i(Z-\alpha Y),i(Z-\alpha Y)].
\end{align*}
The precise identity from that proof is
\[
 \frac{D^2\Omega(Y)[Z,Z]-(p-1)D^2\Omega(Y)[iZ,iZ]}{p\Omega(Y)}
 =\frac{\Delta_\Phi(Z-\alpha Y)}{\Tr(|Y|^r)}
       +(p-1)(r-p)\beta^2.
\]
Both terms on the right are nonnegative. 
For real \(\lambda,\mu\), set
\[
 \delta_s (\lambda,\mu)=
 \begin{cases}
 \dfrac{\lambda|\lambda|^{s-2}-\mu|\mu|^{s-2}}{\lambda-\mu},
      &\lambda\neq\mu,\\[5pt]
 (s-1)|\lambda|^{s-2},&\lambda=\mu.
 \end{cases}
\]
For \(r >2\),
\[
 \delta_r (\lambda,\mu)-(p-1)\delta_r (-\lambda,\mu)\geq 0,
\]
and equality holds precisely when \(\lambda+\mu=0\).
This strictness follows from strict concavity of \(\tanh\) on
\((0,\infty)\).

Consequently, the equality in the proposed Hessian inequality would imply
that the two self-adjoint dilations
\[
 \mathsf A=\begin{pmatrix}0&Y\\Y^*&0\end{pmatrix},\qquad
 \mathsf S=\begin{pmatrix}0&Z-\alpha Y\\Z^*-\overline{\alpha }Y^*&0\end{pmatrix}
\]
anticommute. 
Indeed, in spectral matrix entries of \(\mathsf A\),
the positive defect forces every entry of \(\mathsf S\) with \(\lambda+\mu\neq0\) to vanish. 
The zero--zero block causes no
problem. In particular,
\begin{equation}\label{eq:anticommute}
 Y(Z-\alpha Y)^*+(Z-\alpha Y)Y^*=0,\qquad Y^*(Z-\alpha Y)+(Z-\alpha Y)^*Y=0.
\end{equation}

Take singular vectors \(Yv_i=s_i u_i\), \(s_i>0\).
All these vectors lie in \(\operatorname{Dom}K\). 
To justify this point, \(Y^*Y\) and \(YY^*\) are in \(C^1(K)\); 
every isolated positive-eigenvalue projection $P$ is in \(C^1(K)\) by the resolvent functional calculus. 
Such a projection \(P\) preserves \(\operatorname{Dom}K\). 
Its image of that domain is dense in the finite-dimensional space \(P\mathcal H\), hence equals the whole space. 
Thus \(P\mathcal H\subset\operatorname{Dom}K\).

Put \(a_{ij}=\langle u_i,Ku_j\rangle\) and
\(b_{ij}=\langle v_i,Kv_j\rangle\). Since
\(Z_{ij}=s_j a_{ij}-s_i b_{ij}\), equations
\eqref{eq:anticommute}, for \(s_i\neq s_j\), give
\begin{align*}
 2s_i s_j a_{ij}-(s_i^2+s_j^2)b_{ij}&=0,\\
 (s_i^2+s_j^2)a_{ij}-2s_i s_j b_{ij}&=0.
\end{align*}
The determinant is nonzero, so \(a_{ij}=b_{ij}=0\).
Projecting \eqref{eq:anticommute} to the kernels of \(Y\) and
\(Y^*\) also shows that \(Kv_j\) and \(Ku_j\) have no component
in the corresponding zero singular subspaces. Therefore \(K\)
preserves each nonzero singular-value eigenspace of \(|Y|\).
Each such space is finite dimensional and nonzero. The restriction
of the self-adjoint operator \(K\) to it would have an eigenvector,
contrary to the hypothesis. Equality is impossible.
\end{proof}

\begin{proposition}\label{prop:weylstrict}
Suppose that \(1<p<2\) and assume that \(\theta\) is a nondegenerate real skew-symmetric \(2n\times2n\) matrix.
If $\|U_\theta (f)\|_{L^{p'}(\R_\theta^{2n})} = B_p^{2n} \|f\|_{L^p(\R^{2n})}$, then $f=0$.
If $\|U_\theta^{-1} (X)\|_{L^{p'}(\R^{2n})} = B_p^{2n} \|X\|_{L^{p}(\R_\theta^{2n})}$, then $X=0$.
\end{proposition}

\begin{proof}
Put \(C_p=B_p^{2n}\). 
The localization flow of \cite[Theorem 4.3]{QX} is
\[
 \mathcal L_A(f)=\left(\frac{p A}{\pi}\right)^n
       \int_{\R^{2n}} \norm{U_\theta (f e^{-A|\cdot-s|^2})}_{p'}^p\,ds.
\]
Its initial value is \(\norm{U_\theta f}_{p'}^p\), and it is nondecreasing.
Applying the sharp Weyl inequality to $f e^{-A|\cdot-s|^2}$ gives
\[
 C_p^p\norm f_p^p\leq\mathcal L_A(f)
       \leq C_p^p\norm f_p^p.
\]
Thus $f e^{-A|\cdot-s|^2}$ is an extremizer for almost every \(s\).
Continuity in \(s\) implies $f e^{-A|\cdot-s|^2}$ is an extremizer for every \(s\).

For \(a>0\) and \(b\in\R^{2n}\), set
\[
 Y(a,b)=U_\theta (f e^{-a|\xi|^2+b\cdot\xi}),\quad
 V(a,b)=\norm{Y(a,b)}_{p'}^p,\quad
 I(a,b)=\int |f(\xi)|^p e^{-pa|\xi|^2+pb\cdot\xi}\,d\xi.
\]
The extremality implies \(V(a,b)=C_p^p I(a,b)\). 
Direct differentiation shows
\begin{equation}\label{eq:equalityPDE}
 (\Delta_b+p\partial_a)V(a, b)=0.
\end{equation}
All weighted moments of $f e^{-A|\cdot-s|^2}$ belong to \(L^1\cap L^p\). 
Parameter derivatives therefore exist in \(\mathcal{S}^{p'}\) and in operator norm.

Nondegeneracy of the Weyl relations means that multiplication by \(e^{iv\xi_j}\) is implemented by conjugation with \(e^{ivK_j}\), 
where \(K_j\) is a nonzero linear position or momentum observable. 
Each \(K_j\) has no eigenvalues. 
We have
\[
 Z_j=U_\theta(\xi_j f e^{-a|\xi|^2+b\cdot\xi})=[K_j,Y]
\]
after choosing the implementing generator with the corresponding sign and normalization. 
Moreover \(Y\in C^1(K_j)\).
The same differentiation used in \eqref{eq:heat-subsolution} gives
\[
 (\Delta_b+p\partial_a)V(a,b)
 =\sum_{j=1}^{2n}
 \left(D^2\Omega(Y)[Z_j,Z_j]-(p-1)D^2\Omega(Y)[iZ_j,iZ_j]\right).
\]
The Weyl transform $U_\theta$ is faithful, and
Lemma~\ref{lem:strict} makes every term strictly positive. 
This contradicts \eqref{eq:equalityPDE}. 
The forward inequality has no nonzero extremizer.

The inverse statement follows by norm-attainment duality. 
If \(U_\theta^{-1}(X)\in L^{p'}(\R^{2n})\) attained the inverse norm, choose a
unit \(g\in L^p\) norming \(U_\theta^{-1} (X)\). 
Then
\[
 C_p\norm X_p
 =|\langle U_\theta^{-1}(X),g\rangle|
 =|\tau(X(U_\theta g)^*)|
 \leq\norm X_p\norm{U_\theta g}_{p'}
 <C_p\norm X_p,
\]
a contradiction. 
\end{proof}

\begin{proposition}
Suppose that $\G$ is Kahng's skew-symmetric quantum group and $1< p<2$.
Suppose that $\norm{\F_{p}(X)}_{L^{p'}(\widehat\G)}
 = B_p^{2n+1}\norm X_{L^p(\G)}$.
 Then $X=0$.
\end{proposition}
\begin{proof}
For almost every \(t\neq0\), the source fiber is a nondegenerate
\(2n\)-dimensional Weyl fiber. Proposition~\ref{prop:weylstrict}
therefore makes \eqref{eq:fiberHY} strict whenever \(X_t\neq0\).
If \(X\neq0\), the set of such parameters has positive measure.
Integrating the strict inequality gives
\[
 \norm f_{L^p_tL^{p'}_a}<B_p^{2n}\norm X_p.
\]
The synthesis estimate then implies
\[
 \norm{\F(X)}_{p'}<B_p^{2n+1}\norm X_p.
\]
This proves the proposition.
\end{proof}

\section*{Declaration of generative AI use}

During the preparation of this work, the authors used generative AI
tools to improve the language and readability of the manuscript and
to assist with developing and checking mathematical arguments.
The authors reviewed and revised the AI-assisted content and verified
the mathematical statements and proofs.
The authors take full responsibility for the accuracy and integrity
of the manuscript.

Acknowledgement. The authors are supported by the National Natural Science Foundation
of China (Nos.  and 12601240).
\bibliographystyle{plain}
\bibliography{hyqg}
\end{document}